\documentclass[12pt]{article}

\usepackage{color}

 \usepackage[margin=1in]{geometry}
\usepackage{amsmath,amsthm}
\usepackage{bm}
\usepackage{amsfonts}
\usepackage{amssymb}
\usepackage{color}
\usepackage{graphicx}
\usepackage{natbib}
\usepackage{float}
\usepackage{wrapfig}
\usepackage{hyperref}
\usepackage{siunitx}
\usepackage{caption}
\usepackage{subcaption}
\usepackage{soul}

\newtheorem{theorem}{Theorem}

\newcommand{\ARL}{\mathrm{ARL}}
\newcommand{\EDD}{\mathrm{EDD}}
\newcommand{\E}{\mathbb E}
\newcommand{\Cov}{\mathrm{Cov}}

\title{Sequential change-point detection:\\ Computation versus statistical performance}

\author{
Haoyun Wang$^1$, Yao Xie$^1$\footnote{yao.xie@isye.gatech.edu}\\
        \small $^{1}$School of Industrial and Systems Engineering, Georgia Institute of Technology,\\
        \small Atlanta, Georgia, 30332, USA \\
}
\date{}

\allowdisplaybreaks

\begin{document}

\maketitle

\begin{abstract}

Change-point detection studies the problem of detecting the changes in the underlying distribution of the data stream as soon as possible after the change happens.  Modern large-scale, high-dimensional, and complex streaming data call for computationally  (memory) efficient sequential change-point detection algorithms that are also statistically powerful. This gives rise to a computation versus statistical power trade-off, an aspect less emphasized in the past in classic literature. This tutorial takes this new perspective and reviews several sequential change-point detection procedures, ranging from classic sequential change-point detection algorithms to more recent non-parametric procedures that consider computation, memory efficiency, and model robustness in the algorithm design. Our survey also contains classic performance analysis, which provides useful techniques for analyzing new procedures.

\end{abstract}


\section{Introduction}

Sequential change-point detection has been a classic topic in statistics since the 1920s \citep{shewhart1925application,shewhart1931economic} with motivations in quality control --- the objective is to monitor the manufacturing process by examining its statistical properties and signals a change in quality when the distribution of certain features of the products deviates from the desired one. Since then, change-point detection finds many applications in other fields, including monitoring power networks \citep{chen2015quickest}, internet traffic \citep{lakhina2004diagnosing}, sensor networks \citep{hallac2015network,raghavan2010quickest}, social networks \citep{raginsky2012sequential,peel2015detecting}, medical image processing \citep{malladi2013online}, cybersecurity \citep{tartakovsky2012efficient}, video surveillance \citep{lee2005online}, COVID-19 intervention \citep{dehning2020inferring}, and so on. Recently, there has been much interest in out-of-distribution detection \citep{ren2019likelihood,magesh2022multiple}, which aims to detect the shift of the underlying distribution of the test data from training data, which is related to sequential change-point detection. 


This tutorial considers sequential (also known as the online or quickest, in some context) change point detection problems, where the goal is to monitor a sequence of data or time series for any change and raise an alarm as quickly as possible the change has occurred to prevent any potential loss. The above problem is different from another major type of change-point detection problem, which is offline and aims to detect and localize (possibly multiple) change-points from sequential data in retrospect (see \cite{truong2020selective} for a review). This tutorial focuses on sequential change-point detection. 

\cite{shewhart1925application} introduced a control chart that computes a statistic from every sample, where a low statistic value represents the process still within the desired control. Later it is generalized to compute a statistic for several consecutive samples. \cite{page1954continuous} proposed the CUSUM procedure based on the log-likelihood ratio between the two exactly known distributions, one for the control and one for the anomaly. \cite{moustakides1986optimal,lorden1971procedures} proved that CUSUM has strong optimality properties.  Shiryaev-Roberts (SR) procedure \citep{shiryaev1963optimum,roberts1966comparison} is similar to CUSUM inspired by the Bayesian setting, which is also optimal in several senses \citep{pollak1985optimal,pollak2009optimality,polunchenko2010optimality,tartakovsky2012third}. More recent works on change point detection focus on relaxing the strong assumptions of parametric and known distributions. The generalized likelihood ratio (GLR) procedure \citep{lorden1971procedures,siegmund1995using} aims to detect the change to an unknown distribution assuming a parametric form by searching for the most probable post-change scenario. These classic procedures and their variants can also be found in \cite{basseville1993detection,tartakovsky2014sequential}. \cite{sparks2000cusum} proposed the first adaptive CUSUM, which detects an unknown mean shift by inserting an adaptive estimate of new mean to the CUSUM recursion, followed by \cite{lorden2005nonanticipating,abbasi2019optimal,cao2018sequential,xie2022window}. Such procedures allow an unknown post-change distribution while also having the computational benefit of CUSUM. More recently, \cite{romano2023fast,Poisson-FOCuS2022} present a novel computationally-efficient approach that implements the
GLR test statistic for the Gaussian mean shift and the Poisson parameter shift detection, respectively. The above works belong to parametric change-point detection, i.e., assuming the pre- and post-change distributions belong to the parametric family and detect a certain type of change (for instance, the mean shift and covariance change). There have also been many non-parametric and distribution-free change-point detection algorithms developed, such as kernel-based methods \citep{harchaoui2008kernel,li2019scan,song2022new} and graph-based method \citep{chen2015graph,chu2019asymptotic}.


Sequential change-point detection is different from the traditional fixed-sample hypothesis test. It is fundamentally a repeated likelihood ratio test because of the unknown change point. Moreover, the test statistics are correlated when scanning through time to detect a potential change point. Such correlation is significant, complicates the analysis, and must be explicitly characterized. For instance, in the GLR procedure, we are interested in the probability that the maximum likelihood ratio of all segments of consecutive observations exceeds a certain threshold. In this paper, we review two different analytical techniques to tackle such a challenge: an earlier technique developed based on renewal theory \citep{siegmund1985sequential}, and a more recently developed method based on a change-of-measure technique for the extreme value of Gaussian random fields \citep{yakir2013extremes}.

Despite many methodological and theoretical development for sequential change-point detection, the aspect of (computational and memory) efficiency versus performance tradeoff has been less emphasized. Many change-point detection algorithms are designed, taking such considerations implicitly. This trade-off is becoming more prominent in modern applications due to the need to process large-scale streaming data in real-time and detect changes as quickly as possible. For example, in social network monitoring, one Twitter community will contain at least hundreds of users, let alone the entire network. 

Thus, in this tutorial, we aim to contrast classic and new approaches through the lens of {\it computational efficiency versus statistical performance tradeoff}; we hope this can enlighten developing new techniques for analyzing new algorithms. We highlight the tradeoff by presenting the classic and new performance analysis techniques, ranging from the most standard parametric models to the more recent nonparametric and distribution-free models introduced to achieve computational benefits, such as robustness and flexibility in handling real data. We focus on proof techniques such as renewal theory and the change-of-measure technique, as they are essential for analyzing the classic performance metrics such as the Average Run Length (ARL) and Expected Detection Delay (EDD).

The rest of the paper is organized as follows. Section \ref{sec:setup} describes the formulation of a sequential change-point detection problem and several classic detection procedures. In Section \ref{sec:performance_analysis}, we define the performance metric and provide the classic analysis of CUSUM in renewal theory. Section \ref{sec:adaptive_cusum}, \ref{sec:non_parametric}, \ref{sec:distribution_free} each gives a detailed review of recent literature on managing unknown distributions, from parametric to distribution-free. Section \ref{sec:conclusion} concludes our review.

We would like to acknowledge that there are other recent techniques for analyzing alternative performance metrics, such as \cite{maillard2019sequential,yu2020note,chen2022high}, which we do not discuss in this paper due to space limitations.

\section{Basics of change-point detection and problem setup}\label{sec:setup}


\subsection{Problem setup}

Consider a basic sequential change-point detection problem, where the aim is to detect an unknown change-point. We are given samples $x_1,x_2,\dots$ such that before the change happens, the data distribution follows distribution $f_0 \in \mathcal F_0$, and after the change happens, the distribution shifts to $f_1 \in \mathcal F_1$. At each time step $t$ after collecting the sample $x_t$, we can decide whether to raise an alarm or not. Our goal is to raise an alarm as soon as the change has happened under the false alarm constraint.


The simplest setting is when the sets are singletons, i.e., the pre- and post-change distributions are  $f_0$ and $f_1$, respectively. This happens when we understand the physical nature of the process or there is enough historical data to estimate the pre-change distribution $f_0$. The post-change distribution $f_1$ can be either a target anomaly or the smallest change we wish to detect. 
The anomaly may occur at some time $v\in\mathbb N$, resulting in $x_1,\dots,x_{v-1}\overset{\text{i.i.d.}}{\sim} f_0$ and $x_v,x_{v+1},\dots\overset{\text{i.i.d.}}{\sim} f_1$. If the change never occurs, we write $v = \infty$. Here we assume the change-point $v$ is deterministic but unknown, and the observations before and after the change-point are all independent and identically distributed. 

A change-point detection procedure is a \textit{stopping time} $\tau$, which has the following property. For each $t = 1,2,\dots$, the event $\{\tau>t\}$ is measurable with respect to $\sigma(x_1,\dots,x_t)$, the $\sigma$-field generated by data till time $t$. In other words, whether the procedure decides to stop right after observing the $t$-th sample depends solely on the history, not the future. Common detection procedures are based on the choice of a certain detection statistic, computing the detection statistic using the most recent data to form $T_t$, and stopping the first time that the detection statistics $T_t$ signals a potential change-point (typically by comparing with a predetermined threshold $b$; the choice of the threshold $b$ balances the false alarm and detection delay). Thus, the detection procedure can be defined as a stopping time  
\[\tau = \min\{t\geq 1:T_t>b\}.\] 

Since, in practice, there are usually abundant pre-change samples for us to estimate $f_0$ with high accuracy, here we only consider the case with known pre-change distribution $f_0$. Uncertainties in the pre-change distribution can also be addressed though, and we give such an example in Section \ref{sec:non_parametric} and \ref{sec:distribution_free}. In such cases, it will be harder to control false alarms, and the detection procedure will be more complex computationally in general. We would often like to treat the post-change distribution $f_1$ as unknown since it is due to an unknown anomaly. 
Then instead of considering a single distribution $f_1$, we consider the post-change distribution belonging to a parametric distribution family $f_1(\cdot,\theta)$ with parameter $\theta\in\boldsymbol\Theta_1$. 

\subsection{Computation and robustness considerations versus statistical performance}

Two commonly used performance metrics for change-point detection are the average run length (ARL) and the expected detection delay (EDD), which we describe below and introduce more formally in Section \ref{sec:performance_metrics}. The ARL (related to the false alarm rate) is the expected stopping time under the pre-change distribution (i.e., where there is no change point), and the EDD is the expected stopping time after a change point has happened. Typically, one would choose a threshold to control the ARL to satisfy $\ARL \geq \gamma$ for some chosen large lower bound $\gamma$, and measure the statistical performance by the EDD for a fixed $\gamma$. The well-known lower bound (see, e.g., \cite{lorden1971procedures}) is $\EDD = \log\gamma(1+o(1))/D(f_1\|f_0)$ as $\gamma\to\infty$, where $D(f_1\|f_0)$ is the Kullback-Leibler (KL) divergence between $f_0$ and $f_1$. A larger KL divergence means it is easier to distinguish $f_1$ from $f_0$ and thus a smaller EDD is possible. 

Another way to describe the objective is to minimize
$$
 -\ARL + \lambda~\EDD,
$$
with some hyper-parameter $\lambda$. For some detection procedure $\tau$ which minimizes the above for certain $\lambda$, it also minimizes the EDD among all detection procedures with $\ARL\geq \gamma = \ARL(\tau)$. Most of the existing change-point detection procedures can be regarded as trying to solve the following minimax problem (with slightly varying definitions on the ARL and EDD):
\begin{equation}
\min_{\tau\in \mathcal T}\max_{f_1\in \mathcal F_1} -\ARL_{f_0}(\tau) + \lambda~\EDD_{f_1}(\tau),
\end{equation}
where $\mathcal T$ is the set of stopping times, and $\mathcal F_1$ is the set of possible post-change distributions. We include $f_0,f_1$ in the subscript to show the dependence of the performance metrics on those distributions, but since there is often a large amount of reference data to estimate the pre-change distribution, here we only consider $f_1$ to be unknown. The statistical performance of a detection procedure can then be represented by the value $
\max_{f_1\in \mathcal F_1} -\ARL_{f_0}(\tau) + \lambda~\EDD_{f_1}(\tau).
$
One can expect that as the size of $\mathcal F_1$ grows, a detection procedure either becomes more complex or has worse statistical performance. This leads to the trade-off between computation, model robustness (the size of $\mathcal F_1$), and statistical performance.


Given observations $x_1,x_2,\dots,x_t$, intuitively, to utilize observations and achieve the best statistical performance fully, we would use all the past samples in the detection statistic -- as is in the case for the CUSUM and GLR procedures. However, this can become prohibitive in practice as $t$, the duration we have run the detection procedure grows larger and larger (unless the algorithm is fully recursive such as CUSUM). 

A practical online change-point detection algorithm should have constant computation complexity  and memory requirement $O(1)$ per iteration (unit time), but this is not likely to be statistically powerful for various cases. However, in some situations, the statistically powerful algorithm will require $O(t)$ computation per iteration, which grows with time and thus is not practical. Therefore, we constantly face a computation and statistical performance tradeoff. Due to this consideration, a commonly adopted simple (yet effective in many cases) strategy is to use {\it sliding window}: which stores historical data within a sliding window of length $w$ and computes the statistic $T_t$ using data in the sliding window $(x_{t-w+1},x_{t-w+2},\dots,x_{t})$ for every $t\geq w$; this way, both memory and computational complexities are constant in the duration $t$. We wish to have $w$ as small as possible to minimize the algorithm's memory complexity and computation complexity. However, $w$ also cannot be too small to sacrifice performance. Thus, the critical question is to decide the window length $w$--how much data needs to be remembered. We will see a trade-off between the memory complexity with respect to the window length $w$ and statistical performance in Section \ref{sec:other_procedures}.

The popularity of CUSUM in practice is probably due to that it achieves constant memory and computational complexity (only use the current sample). Furthermore, its statistical performance is asymptotically optimal (which we will specify more precisely later). However, this optimality of CUSUM requires precise knowledge of the pre-change and post-change distributions, which is not robust to model misspecification. In improving the model's robustness by relaxing such requirements (especially the assumption of known post-change distribution), the tradeoff arises between computation and memory complexity versus statistical performance. For example in window-limited GLR, the memory requirement is $O(w)$ where the window length is directly related to the gap between the pre-change distribution $f_0$ and $\mathcal F_1$ in order to have good statistical performance. And the computation complexity, in general, is at least $O(w^2)$ per time unit.

Another direction in developing robust change-point detection procedures is to utilize distribution-free methods and non-parametric statistics. Arguably, when the distributional models can be specified more or less reasonably, the non-parametric models are not needed, and they may not be as good as the parametric change-point detection algorithms. However, they gain robustness when the data distributions are not easy to specify using parametric models. While the construction of the non-parametric detection statistics can be straightforward in many cases (such as the kernel-based sequential change-point detection procedure \cite{li2019scan}), the performance analysis is much harder than the parametric cases due to a lack of handle through the probability density function in parametric distributional models; the distribution of the non-parametric statistics can be unknown functional form. The asymptotic optimality is also harder to analyze in the sense that a meaningful lower bound is unclear. There are several recent interesting works in this area, and we discuss them in Section \ref{sec:non_parametric} and Section \ref{sec:distribution_free}. 

\begin{figure}[htbp]
\centering\includegraphics[width=0.8\textwidth]{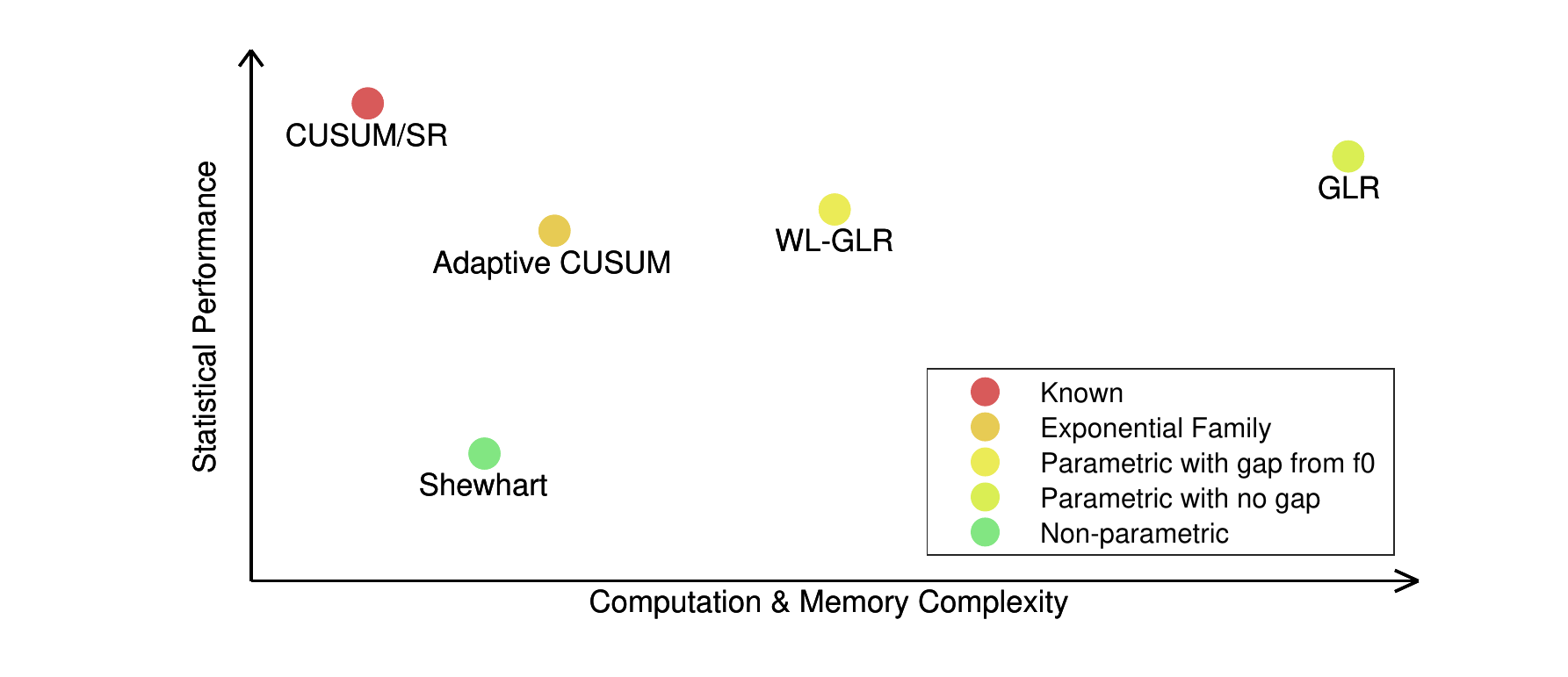}
\caption{Illustration of the statistical performance, computational and memory efficiency, model robustness trade-off for the classic procedures in sequential change detection. The color shows the assumption on post-change distribution.}
\label{fig:trade-off-comparison}
\end{figure}

In the following, we will describe several common procedures. The comparisons of the detecting statistics are summarized in Figure \ref{fig:trade-off-comparison} and Table \ref{table:pros_cons}.

\subsection{Classic CUSUM procedure}

The CUSUM procedure is derived based on likelihood ratios. For an assumed change-point location $\nu$, the log-likelihood ratio for the hypothesis 
$H_0: x_1,x_2,\dots\sim f_0$ versus the alternative hypothesis $H_v: x_1,\dots,x_{v-1}\sim f_0, x_v,x_{v+1}\dots\sim f_1$ is given by
\[
\sum_{i = v}^t \log\left(\frac{f_1(x_i)}{f_0(x_i)}\right) 
\]
Since change-point $\nu$ is unknown, the detection statistic needs to consider the maximization of the above with respect to all possible change-point locations. 
This gives rise to the following detection statistic for each time $t$
\begin{equation}
T_t^{\rm CUSUM} = \max_{1\leq k\leq t+1} \sum_{i=k}^t \log\left(\frac{f_1(x_i)}{f_0(x_i)}\right).
\label{eq:cusum_stat_def}
\end{equation}
Here when $k=t+1$, we set the empty summation to 0. The CUSUM procedure computes the detection statistic for each time $t$, stops the first time that the detection statistic exceeds a certain threshold $b > 0$, and claims there has been a change in the past:
\begin{equation}
\tau^{\rm CUSUM}(b) = \inf\{t: T_t^{\rm CUSUM} > b\}.
\label{eq:cusum_def}
\end{equation}

The popularity of CUSUM is possibly due to the following recursive computation of the detection statistic. 
Let $z_i = \log(f_1(x_i)/f_0(x_i))$ be the increment in the log-likelihood ratio, $i=1,2,\dots,$ and $ S_t = \sum_{i=1}^tz_i$, $\forall t\geq 0$. Then the CUSUM statistic has the following recursive expression because the partial sums starting from each potential change-point $k$ shares the same increment when updating from time step $t$ to $t+1$:
\begin{align}
T_{t+1}^{\rm CUSUM}  = \max_{0\leq k\leq t+1}( S_{t+1} -  S_k) 
=&\ \max\left\{z_{t+1} + \max_{0\leq k\leq t} (  S_t  -  S_k),0\right\}\nonumber \\
=&\ \max\left\{z_{t+1} + T_{t}^{\rm CUSUM},0\right\},
\label{eq:simple_recursion}
\end{align}
with $T_0^{\rm CUSUM} = 0$. This convenient recursive CUSUM evaluation (\ref{eq:simple_recursion}) means that we do not need to remember any data and merely update the detection statistic every time we observe a new sample. More precisely, for each time slot, one only needs to compute the new log-likelihood ratio increment $z_i$ and update the CUSUM statistic $T_{t}^{\rm CUSUM}$ according to \eqref{eq:simple_recursion}. In other words, both the computation and memory complexity per update are constant for CUSUM. 

Moreover, CUSUM's statistical performance enjoys asymptotic optimality \citep{lorden1971procedures}, which we will explain in Section \ref{sec:performance_metrics}. It is exactly optimal shown later by \citep{moustakides1986optimal}. Largely speaking, for CUSUM to enjoy a good performance, we need the property that the expected value of the increment term before the change happens is negative ($\mathbb E_0[z_t] < 0$), and positive after the change ($\mathbb E_1 [z_t] > 0$). This can be verified using Jensen's inequality (when the data distribution specifications are precise). 

However, a known drawback of CUSUM is that it requires both the pre-change and post-change distributions to be parametric and specified exactly, which can be too restrictive for real-world applications. When the true distribution deviates from the assumed distributions, CUSUM is no longer optimal and suffers from performance loss depending on the level of model misspecification. For modern data, especially high-dimensional data, having an exact specification of data distribution is difficult, and performance degradation of CUSUM due to model mismatch becomes inevitable. 

Multiple variants have been developed to make CUSUM more robust to unknown post-change distributions. One possibility is to run multiple CUSUM procedures in parallel, each detecting against a different possible anomaly outcome \citep{lorden1973detection,lucas1982combined,tartakovsky2005asymptotic,xie2013sequential}. Another approach is to consider the so-called {\it adaptive CUSUM}, to adaptively estimate the post-change distribution while running the CUSUM recursion \citep{sparks2000cusum,lorden2005nonanticipating,abbasi2019optimal, xie2020sequential}. Uncertainties in the pre-change distribution can be treated similarly \citep{pollak1991sequential,krieger1999detecting,mei2006sequential}; however, to the best of our knowledge, there is yet a study regarding how to obtain recursive expression when the pre-change distribution is unknown. A more detailed review of adaptive CUSUM and its recent progress can be found in Section \ref{sec:adaptive_cusum}.

\subsection{Other classic procedures}

\label{sec:other_procedures}
\textbf{Shiryaev-Roberts procedure.} The Shiryaev-Roberts (SR) procedure, similar to the CUSUM, is based on the log-likelihood ratio between the completely specified pre-change and post-change distribution. It is first inspired by putting a prior on the change-point $v$ and assuming it follows a geometric distribution \citep{shiryaev1963optimum}. The SR procedure can now be seen as the resulting procedure when the geometric mean goes to infinity. This is reflected in the detection statistic, where the SR procedure takes the summation instead of the maximum over all the potential change-point $k$ (as done in CUSUM),
$$
T_t^{\rm SR} = \sum_{k=1}^{t} \prod_{i=k}^t\frac{f_1(x_i)}{f_0(x_i)}.
$$
And the recursion (similar to CUSUM) is in the form of
$$
T_{t+1}^{\rm SR} = (1+T_t^{\rm SR})\frac{f_1(x_{t+1})}{f_0(x_t)}.
$$
The SR procedure is exactly optimal under the integral expected detection delay \citep{pollak2009optimality} and later \cite{tartakovsky2012third} proved that the SR procedure is third order asymptotically optimal following \cite{pollak1985optimal}'s definition of the EDD. They also compared different initial values of $T_0^{\rm SR}$ and discussed the detection delay conditioned on the change-point $k$. See \cite{polunchenko2012state} for a review on sequential change-point detection with known pre- and post-change distributions.

\vspace{.1in}
\noindent
\textbf{Generalized Likelihood Ratio (GLR) procedure.} The GLR procedure is adopted when the post-change distribution is unknown but still parametric. The test statistic scan through the log-likelihood ratio over all the potential change-point $k$ and for each assumed $k$ and maximum likelihood estimator of the post-change parameter is used in forming the likelihood ratio: 
$$
T_t^{\rm GLR} =\max_{1\leq k\leq t} \sup_{\theta\in\boldsymbol\Theta_1} \sum_{i=k}^t\log\left(\frac{f_1(x_i,\theta)}{f_0(x_i)}\right).
$$
The GLR procedure is asymptotically optimal in its statistical performance and handles the problem of unknown post-change distribution automatically by searching through all possible parameters in the feasible region. \cite{lorden1971procedures} proves the first-order asymptotic optimality for the GLR procedure with univariate exponential family, and later \cite{siegmund1995using} provides a more precise characterization of the ARL on the problem of a mean shift with known variance, which turns out to be very useful in analyzing modern change-point detection procedures as well. We will review this method in Section \ref{sec:distribution_free}. Unfortunately, in general, $T_{t+1}^{\rm GLR}$ cannot be updated recursively from $T_t^{\rm GLR}$ and the computation and memory needed per update is often at least linear in $t$. An exception with univariate exponential family can be found in \cite{romano2023fast} where the complexity is reduced to $O(\log t)$. 

\vspace{.1in}
\noindent
\textbf{Window-limited GLR.} A more computationally efficient procedure than the above vanilla GLR is the window-limited GLR \citep{willsky1976generalized} developed by taking the maximum over all potential change-points within a sliding window of fixed length $w$: 
$$
T_{t}^{\rm WL-GLR} = \max_{t-w\leq k\leq t}\sup_{\theta\in\boldsymbol\Theta_1} \sum_{i=k}^t\log\left(\frac{f_1(x_i,\theta)}{f_0(x_i)}\right).
$$
At each time step, the computational and memory requirement of window-limited GLR is constant with respect to $t$. It appears that by forgetting previous samples before time $t-w$, we may lose information. But it can be shown that with a proper window length $w$,  the window-limited GLR can still be asymptotically optimal \citep{lai1999efficient}; this requirement is $w\geq \log \gamma/D_{\min}$, where $\gamma$ is the ARL and $D_{\min}$ is the smallest KL-divergence between $f_0$ and potential $f_1$ we want to detect. Because the window-limited GLR still scans through all potential change points within the sliding window, the computational cost scales with $w$; in general, it is reasonable to assume at least $O(w)$ operations are needed to find the supremum over $\theta\in\boldsymbol \Theta_1$. There are $w$ potential change points per update so the computation complexity would be at least $O(w^2)$. When working with exponential family, however, the computation complexity is reduced to $O(w)$ because the partial sum $\sum_{i=k}^t x_i$ is a sufficient statistic for computing the maximum likelihood ratio for each potential change point $k$.

\vspace{.1in}
\noindent
\textbf{Shewhart chart} is one of the earliest sequential change-point procedures \citep{shewhart1925application,shewhart1931economic}. Still, many recent detecting procedures fall into this type due to its simplicity: an offline test can be converted into a sequential change-point detection procedure easily by applying the test to a {\it sliding window of samples}. More precisely, the detection statistic $T_t$ using a sliding window of fixed length $w$. For instance, the simplest example is by letting $T_t$ be the average of the past $w$ samples to detect a mean shift. If the problem is parametric, it can be the generalized likelihood ratio: 
\begin{equation}
T_t^{\rm SH-GLR} = \sup_{\theta\in\boldsymbol\Theta_1} \sum_{i=t-w}^t\log\left(\frac{f_1(x_i,\theta)}{f_0(x_i)}\right),
\label{eq:sh-glr}
\end{equation}
And again when dealing with exponential families, the above turns into an explicit function of the partial sum $\sum_{i=t-w}^t x_i$. Another example of Shewhart chart is based on the score statistic (e.g., \cite{chen2017textsf}) utilizing the locally most powerful score statistic, 
$$
T_{t}^{\rm Score} = \frac{1}{2w}\nabla_\theta^T \left(\sum_{i=t-w}^t \log f_1(x_i,\theta)\right)  \mathcal I^{-1}(0)
\nabla_\theta \left(\sum_{i=t-w}^t \log f_1(x_i,\theta)\right),
$$
where $\mathcal I(0)$ is the Fisher information \cite{casella2021statistical} at $\theta = 0$ and it can be pre-computed.

The score statistic can sometimes lead to simple detection statistics. Since we only consider the gradient of the likelihood ratio at the single parameter value under the null hypothesis, sometimes this avoids solving optimization problems or matrix inversion. Similar to the window-limited GLR, the choice of the window length $w$ depends on the smallest change we are interested in detecting.  To further simplify computation, one can replace computing the statistic $T_t$ at every time step with updating it every $\delta$ time unit for $\delta$ up to $w$. 

We would like to emphasize that a distinction is that the Shewhart detection statistic, unlike CUSUM, does not scan over all potential change-point locations. Because the Shewhart chart no longer scans over potential change points, its statistical performance is not asymptotically optimal with respect to the ARL EDD metric (with the exception that if we consider a different metric: maximizing the probability of detection \cite{Pollak2013,moustakides2014multiple}). It is commonly believed that CUSUM is more statistically powerful for detecting small changes.


\begin{table}[htbp]
    \centerline{
    \small
    \begin{tabular}{|c|c|c|c|c|}
        \hline
        & Model Robustness & \multicolumn{2}{c|}{Computation \& Memory Efficiency} & Asymptotic optimality\\\hline
        CUSUM/SR & $f_1$ is known & $O(1)$ & $O(1)$ & $\checkmark$\\
        Adaptive CUSUM/SR & Exponential family & $O(1)$ & $O(1)$ & $\checkmark$\\
        GLR & Parametric& $\geq O(t)$ & $O(t)$ & $\checkmark$\\
        Window-limited GLR & Parametric & $\geq O(w)$ & $O(w)$ & $\checkmark$\\
        Shewhart chart & Can be nonparametric & Varies & $O(w)$ & $\times$\\\hline
    \end{tabular}
    }
    \vspace{.1in}
    \caption{Comparisons of classic parametric change-point detection procedures. Model robust means whether the detection procedure requires an accurate post-change distribution as an input. Efficiency includes both computational and memory requirements per time unit, and we say the procedure is efficient if both are a constant with respect to the duration $t$ we monitor the process. Asymptotic optimality is reached when the performance satisfies the theoretical lower bound given in Section \ref{sec:performance_metrics} Theorem \ref{thm:theoretical_lower_bound}. Note that the Shewhart chart contains a broad range of detecting procedures and can be based on either parametric or non-parametric statistics. 
    }
    \label{table:pros_cons}
\end{table}

\section{Statistical performance}
\label{sec:performance_analysis}

We start by reviewing statistical performance analysis for change-point detection algorithms. The techniques used for statistical performance analysis include renewal theory, Wald's identity, and change-of-measure techniques  (see, e.g., \cite{page1954continuous,lorden1971procedures,lorden1973detection,pollak1985optimal,siegmund1985sequential,lai1998information}; a recent survey in quickest change-point detection can be found in \citep{xie2021sequential}).  We review such techniques because they are still useful for analyzing new change-point detection procedures \citep{tartakovsky2019sequential}.

\subsection{Standard performance metrics: ARL and EDD}
\label{sec:performance_metrics}

The performance of a change-point detection procedure with stopping time $\tau$ can be evaluated by the average stopping time after a change with a false alarm constraint. Let $\mathbb P_v,\mathbb E_v$ be the probability measure and expectation given time of change $v=1,\dots$ (with $\nu = \infty$ denoting no change). 

Consider the following problem. We would like to find a detection procedure that can minimize the worst-case expected detection delay \citep{pollak1985optimal}:
\begin{equation}
\EDD(\tau) := \sup_{v\geq 1} \mathbb E_v[\tau - v + 1|\tau\geq v],
\label{eq:def_edd}
\end{equation}
subject to average run length
\begin{equation}
\ARL(\tau) := \mathbb E_0[\tau]\geq \gamma
\label{eq:def_arl}
\end{equation}
for some large constant $\gamma > 0$. (Note that there are other definitions of the worst-case detection delay, such as that in \cite{lorden1971procedures}, which we do not consider here.) The theoretical lower-bound of the EDD has related to the Kullback-Leibler (KL) divergence between the post-change and pre-change distributions:
$$
D(f_1\|f_0) = \E_1[\log(f_1(x)/f_0(x))].
$$

\begin{theorem}[Lower bound for EDD \citep{pollak1985optimal,lai1998information}]
$$
\inf_{\tau: \ARL(\tau)\geq \gamma} \EDD(\tau) \geq \frac{\log \gamma}{D(f_1\|f_0)}(1+o(1)),
$$
as $\gamma\to\infty$.
\label{thm:theoretical_lower_bound}
\end{theorem}

\noindent It can be shown that the CUSUM procedure satisfies
\begin{equation}
    \ARL(\tau^{\rm CUSUM}(b)) = \Theta(e^b),\quad \EDD(\tau^{\rm CUSUM}(b)) = \frac{b}{D(f_1\|f_0)}(1+o(1)),
    \label{eq:arl_edd}
\end{equation}
as the threshold $b\to \infty$,  and hence it is asymptotically optimal; here, $\Theta$ is the big-theta notation. In the next subsection, we provide the analysis leading to \eqref{eq:arl_edd} for the basic i.i.d. setup. For non-i.i.d. data, \eqref{eq:arl_edd} still holds when the average log-likelihood ratio between $\mathbb P_0$ and $\mathbb P_v$ converges to some constant in probability as $t\to\infty$ \citep{lai1998information}.

\subsection{Analysis via renewal theory and Wald's identity}
\label{subsec:arl_edd_renewal_theory}


Now we demonstrate the analysis of the ARL and EDD approximation \eqref{eq:arl_edd} for the CUSUM procedure when the data are i.i.d. By the recursive rule \eqref{eq:simple_recursion}, $(T_t^{\rm CUSUM})_{t\geq 1}$ in CUSUM can be viewed as a random walk. Each time the random walk is increased by $\log(f_1(x_{t})/f_0(x_{t}))$, with the exception that if the path steps below 0, there is a renewal and the path is reset to 0. The random walk ends the first time when the path exceeds the threshold $b$.

To analyze the EDD of CUSUM procedure $\tau^{\rm CUSUM}(b)$ for i.i.d. data, first, we have the following fact that greatly simplifies EDD's computation. For any $b\geq 0$, for CUSUM procedure $\tau^{\rm CUSUM}(b)$ (which is abbreviated as $\tau$ in the informal proof to simplify notation)
\begin{equation}
\sup_{v\geq 1} \mathbb E_v[\tau - v + 1|\tau\geq v]
= \mathbb E_1[\tau].
\label{eq:edd_at_1}
\end{equation}
\begin{proof} (Informal)
We first show that for any $v\geq 1$, 
$$
\mathbb E_v[\tau - v +1|\tau\geq v]\leq \mathbb E_1[\tau].
$$
Thus, the conditional expectation of detection delay in EDD's definition \eqref{eq:def_edd} reaches its supremum at $v = 1$. Conditioned on $\tau\geq v$ and $T_{v-1}^{\rm CUSUM}$, $\tau$ is decided by the random walk starting from $T_{v-1}^{\rm CUSUM}\geq 0$ with increments $z_v,z_{v+1},z_{v+2},\dots,$ which are function values of i.i.d. random variables $x_v,x_{v+1},\dots,\sim f_1$. By the independence between samples, the increments $z_v,z_{v+1},\dots,$ are independent of $T_{v-1}^{\rm CUSUM}$ and the event $\tau\geq v$. For every $v\geq 2$ and $s\geq 0$, there must be $$
\E_v\left[\tau - v + 1\big|\tau\geq v,T_{v-1}^{\rm CUSUM} = s\right]\leq \E_v\left[\tau - v + 1\big|\tau\geq v,T_{v-1}^{\rm CUSUM} = 0\right]
$$
because the random walk with renewal must be at least as large as if it starts from $T_{v-1}^{\rm CUSUM} = 0$, and hence stops earlier. Also note that the random walk $(T_{v-1+t}^{\rm CUSUM})_{t\geq 0}$ conditioned on $T_{v-1}^{\rm CUSUM} = 0$ under $\mathbb P_v$ is identically distributed with $(T_t^{\rm CUSUM})_{t\geq 0}$ starting from $T_0^{\rm CUSUM} = 0$ under $\mathbb P_1$, so their expected stopping time should be the same, which gives 
\[
\E_v\left[\tau - v + 1|\tau\geq v,T_{v-1}^{\rm CUSUM} = 0\right] = \E_1[\tau].
\]
Then
\begin{align*}
\E_v\left[\tau - v + 1\big|\tau \geq v\right] = &\ \E_v\left[\E_v\left[\tau - v + 1\big|\tau\geq v,T_{v-1}^{\rm CUSUM}\right]\big|\tau\geq v\right]\\
\leq&\  \E_v\left[\tau - v + 1\big|\tau\geq v,T_{v-1}^{\rm CUSUM} = 0\right]\\
=&\ \E_1[\tau].
\end{align*}
\end{proof}

Now with \eqref{eq:edd_at_1}, we can compute the EDD using renewal properties under $\mathbb P_1$. Let $N_1$ be the first time the random walk $(T_t^{\rm CUSUM})_{t\geq 1}$ falls out of the interval $(0,b)$. If $T_{N_1}^{\rm CUSUM}\geq b$, the random walk stops, otherwise there is a renewal, and $T_{N_1}^{\rm CUSUM}$ is reset to be 0. We can continue and define $N_2,N_3,\dots$ to be the time the renewed random walk $(T_{N_1 + t}^{\rm CUSUM})_{t\geq 1}$ first falls out of $(0,b)$, until finally for some $M$, $T_{N_1+\dots+N_M}^{\rm CUSUM}\geq b$. Since the increments $z_t,t=1,2,\dots$ are i.i.d., so are the random walks after renewal. This means $N_1,N_2,\dots$ are i.i.d, and the conditional probability $\mathbb P_1(M = m|M\geq m)$ is a constant with respect to $m$. So $M$ follows a geometric distribution with expectation $1/\mathbb P_1(M=1)$. Then we need Wald's identity on the expectation of summation of random variables with a stopping time.

\begin{theorem}[Wald's identity]
For i.i.d. random variables $X_1,X_2,\cdots$ and a stopping time $\tau$, if $\E[X_1],\E[\tau]<\infty$, then
$$
\E\left[\sum_{i=1}^\tau X_i\right] = \E[X_1]\E[\tau].
$$
\end{theorem}

\noindent By the Wald's identity,
\begin{equation}
\E_1[\tau^{\rm CUSUM}(b)] = \E_1\left[\sum_{i=1}^MN_i\right] = \E_1[N_1]\E_1[M] = \E_1[N_1]/\mathbb P_1(M = 1).
\label{eq:random_walk_edd}
\end{equation}
For the ARL, the same argument applies with the probability measure $\mathbb P_0$ and 
\begin{equation}
\E_0[\tau^{\rm CUSUM}(b)] = \E_0[N_1]/\mathbb P_0(M=1).
\label{eq:performance_random_walk}
\end{equation}
Now we have established that the ARL equals to $\E_0[N_1]/\mathbb P_0(M=1)$ and the EDD equals to $\E_1[N_1]/\mathbb P_1(M=1)$. Next, we analyze these values using martingale properties. Let $(S_t)_{t\geq 0}$ be the random walk with increments $z_1,z_2,\dots$ without renewal, i.e. $S_t = \sum_{i=1}^t z_i,~\forall t\geq 0$. Let the stopping time
$$
N = N_1 = \min\{t\geq 1: S_t\notin (0,b)\}.
$$
\textbf{Ladder variables.} For now, the increments $z_i$ can be either positive or negative. To connect with the renewal theorem, we introduce the ladder variable $z^+$, which is the value when $(S_t)_{t\geq 0}$ first hits the positive axis, and let 
$$\tau_+ = \inf\{t:S_t>0\}.$$
By the law of large numbers, since $\E_1[z_1]>0$, $\tau_+$ must be finite under $\mathbb P_1$. So $z^+  = S_{\tau_+}$ is well defined under $\mathbb P_1$. We also define two other stopping times for later use. $$\tau_b = \inf\{t:S_t\geq b\}, \quad \tau_- = \inf\{t:S_t< 0\}.$$
As we'll see later, the value of interest is the overshoot $S_{\tau_b}-b$ under $\mathbb P_1$, and it suffices to consider only the distribution of $z^+$. This is because $(S_t)_{0\leq t\leq \tau_b}$ can break down into several pieces, each ending when the random walk reaches a new highest. Therefore, $S_{\tau_b}$ can be expressed as the sum of i.i.d. copies of $z^+$. 

\noindent
\textbf{Approximating ARL.}
For $i = 0$, $(\exp(S_t))_{t=0}^{N}$ is a martingale under $\mathbb P_0$ because for every $0\leq t< N$, 
$$
\E_0[\exp(S_{t+1})|\exp(S_t)] = \exp(S_t)\E_0[f_1(x_{t+1})/f_0(x_{t+1})] = \exp(S_t).
$$
By Doob's Martingale inequality,
\begin{equation}
\mathbb P_0(S_{N}\geq b) = \mathbb P_0(\exp(S_{N})\geq e^b)\leq \frac{\E_0[\exp(S_0)]}{e^b} = e^{-b}.
\label{eq:no_overshoot}
\end{equation}
And obviously $N\geq 1$, so
\begin{equation}
\ARL(\tau^{\rm CUSUM}(b)) = \E_0[N]/\mathbb P_0(S_{N}\geq b) \geq  e^b.
\label{eq:ARL_lower_bound}
\end{equation}
To better approximate the ARL and improve the multiplicative constant, first, we look at $\E_0[N]$. As $b\to\infty$, $N\to\tau_-$ and $\E_0[N]\to \E_0[\tau_-]$ because of the monotone convergence theorem. The divider in \eqref{eq:performance_random_walk} can be evaluated using a change-of-measure trick,
\begin{align}
\mathbb P_0(S_{N}\geq b) = \E_0[\boldsymbol 1\{S_{N}\geq b\}] = &\ \E_1[\exp(-S_{N})\boldsymbol 1\{S_{N}\geq b\}] \nonumber\\
=&\  e^{-b}\E_1[\exp(-(S_{N} - b))\boldsymbol 1\{S_{N}\geq b\}],
\label{eq:change_of_measure}
\end{align}
where $\boldsymbol 1\{\cdot\}$ is the indicator function of an event. The second equality is due to the fact that $\exp(-S_{N})$ is the likelihood ratio between $\mathbb P_0$ and $\mathbb P_1$. We use the renewal theorem to estimate the right-hand side of \eqref{eq:change_of_measure}.
\begin{theorem}[Renewal theorem] \cite[Chapter 8]{siegmund1985sequential}
For any $y\geq 0$, 
$$
\lim_{b\to\infty}\mathbb P_1(S_{\tau_b}- b > y)= \E_1[z^+]^{-1}\int_{y}^\infty\mathbb P_1(z^+ > x)dx.
$$
\label{thm:renewal}
\end{theorem}
This result specifies the asymptotic tail distribution of the ``over-shoot'' term as $b\to\infty$, useful in deriving many properties of the stopping time. Then the right hand side of \eqref{eq:change_of_measure} can be approximated using
\begin{align*}
\E_1[\exp(-(S_{N}-b))\boldsymbol 1\{S_{N}\geq b\}] =&\ \mathbb P_1(S_{N}\geq b) \E_1[\exp(-(S_{N}-b))|S_{N}\geq b]\\
\approx &\ \mathbb P_1(S_{N}\geq b)\E_1[\exp(-(S_{\tau_b}-b))]
\end{align*}
as done in \cite{lorden1973detection}. When $b\to\infty$, $\mathbb P_1(S_{N}\geq b)\to \mathbb P_1(\tau_- = \infty)$, and 
$$
\E_1[\exp(-(S_{\tau_b}-b))]\to \E_1[z^+]^{-1}(1-\E_1[\exp(-z^+)])
$$ as a result of the renewal theorem above. So putting the results above together, we have
$$
\ARL(\tau^{\rm CUSUM}(b)) \approx \frac{\E_0[\tau_-]\E_1[z^+]}{\mathbb P_1(\tau_- = \infty)(1-\E_1[\exp(-z^+)])}e^{b}.
$$


\noindent
\textbf{Approximating EDD.} Since $(S_t - D(f_1\|f_0)t)_{t\geq 0}$ is a martingale under $\mathbb P_1$, by the optional stopping theorem,
$$
\E_1 [S_{\tau_b} - D(f_1\|f_0)\tau_b] = 0,
$$
and $\E_1[\tau_b] = (b+\E_1[S_{\tau^{\rm CUSUM}}-b])/D(f_1\|f_0)$. Note that $S_t \geq T$ for any $t$, so CUSUM must stop no later than $\tau_b$, 
$$
\EDD(\tau^{\rm CUSUM}(b))\leq \E_1[\tau_b]=\frac{b+\E_1[S_{\tau_b}-b]}{D(f_1\|f_0)},
$$
where as a result of Theorem \ref{thm:renewal},
\begin{align*}
\lim_{b\to\infty}\E_1[S_{\tau_b}-b] =&\ \lim_{b\to\infty}\int_0^\infty \mathbb P_1(S_{\tau_b} -b > y) dy = \int_0^\infty \lim_{b\to\infty}\mathbb P_1(S_{\tau_b} -b > y) dy\\
=&\ \E_1[z^+]^{-1} \int_0^\infty\int_y^\infty \mathbb P_1(z^+>x)dxdy = \E_1[z^+]^{-1} \int_0^\infty x\mathbb P_1(z^+>x)dx\\
=&\ \E_1[(z^{+})^2]/2\E_1[z^+].
\end{align*}
Here we assume we can change the order of taking limit and integral but this is not an issue using a more complete version of the renewal theorem. For details, we refer to \cite{siegmund1985sequential}.

For the stopping time $N_1$, similarly, there is
$
\E_1[S_{N} - D(f_1\|f_0)N] = 0,
$ and
$$
\E_1[N] = \frac{\E_1[S_{N}]}{ D(f_1\|f_0)}\geq \frac{\E_1[S_{N}\mathbf 1\{S_{N}\geq b\}]}{D(f_1\|f_0)}\geq \frac{\mathbb P_1(S_{N}\geq b)b}{D(f_1\|f_0)},
$$
with \eqref{eq:random_walk_edd}, together there is
$$
\frac{b}{D(f_1\|f_0)}\leq \EDD(\tau^{\rm CUSUM}(b)) \leq  \frac{b+\E_1[(z^{+})^2]/2\E_1[z^+]+o(1)}{D(f_1\|f_0)}.
$$

\section{Adaptive CUSUM}
\label{sec:adaptive_cusum}
When the post-change distribution belongs to some family parametrized by $\theta\in \boldsymbol\Theta_1$, adaptive CUSUM provides another framework other than GLR.
At each time slot, the statistic $T_t$ is computed using the CUSUM recursive rule \eqref{eq:simple_recursion}, where the increment $z_t$ is the log-likelihood ratio $\log(f_{1}(x_t,\hat\theta)/f_0(x_t))$ for some estimator $\hat\theta_t$. Such method first appears in sequential hypothesis testing where the alternative hypothesis is composite, and a non-anticipating estimator replaces the simple alternative distribution in the sequential likelihood ratio test \citep{robbins1972class,robbins1974expected,pavlov1991sequential}. \cite{tartakovsky2014nearly} compared the adaptive likelihood ratio with the generalized likelihood ratio type-of method, and found the former more convenient in controlling the probability of false alarm. In change-point detection, the key difference between adaptive CUSUM and (even Shewhart-chart type of) GLR is that 
that while both searches for the most probable post-change parameter $\hat\theta_t$, GLR estimates the post-change parameter and finds evidence for the change in the same window, and adaptive CUSUM do these separately. It turns out that we need fewer observations to estimate the post-change parameters than to find evidence, and this can be one reason why adaptive CUSUM outperforms GLR in computational aspects.

An earlier work we found that uses adaptive CUSUM is  \citep{sparks2000cusum}, which considers the case where the change-point causes a mean shift with an unknown scale. They proposed to estimate the shifted mean $\hat\theta_t$ using the exponentially weighted moving average (EWMA) and update the CUSUM statistic using the new sample $x_{t}$ and $\hat\theta_t$ as the targeted change. Then the estimated change is updated using $\hat\theta_{t+1} = \alpha x_t + (1-\alpha)\hat\theta_{t}$ for some $\alpha\in(0,1)$. Following the same spirit, \cite{cao2018sequential} propose to update the estimate $\hat\theta_t$ using online convex optimization algorithms to manage exponential families while the estimated change $\hat\theta_t$ remains to be computed recursively. More recently, \cite{xie2022window} studied a window-limited CUSUM which can be applied to general distribution families of parametric forms. This procedure estimates the $\hat\theta_t$ in a sliding window of length $w$, and the statistic is given by the recursion
\begin{equation}
T_{t+1}^{\rm WL-CUSUM} = \max\left\{T_{t}^{\rm WL-CUSUM} + \log\frac{f_1(x_{t+1},\hat\theta_{t+1}^w)}{f_0(x_{t+1})},0\right\},
\label{eq:wl-cusum_recursion}
\end{equation}
where $\hat\theta_{t}^w$ is a consistent estimator (e.g., the maximum likelihood estimator) of $\theta$ using a sliding window of past $w$ samples: $x_{t-w},\dots,x_{t-1}$. Unlike the window-limited GLR where $w \sim \Theta (\log \gamma)$ needs to be in the same order with EDD (or the threshold $b$) to collect enough information to detect a change, the window-limited CUSUM has smaller memory requirement $w \sim \Theta (\sqrt{ \log \gamma})$, possibly explained by the fact that in window-limited CUSUM we only need enough samples to obtain a reasonably good estimate of $\hat\theta_t$. The optimal window length is explained by the following upper bound on the EDD, and with this window length window-limited CUSUM is asymptotically optimal in its statistical performance:
%
\begin{theorem}[EDD of window-limited CUSUM](\cite{xie2022window})
$$
\EDD \leq \frac{b + \hat J_0/\hat D(f_{1,\theta}\|f_0) + (b\hat J_0/\hat D(f_{1,\theta}\|f_0) )^{1/2} + w D(f_{1,\theta}\|f_0) + (D(f_1\|f_0)w\hat J_0/\hat D(f_{1,\theta}\|f_0) )^{1/2}}{\hat D(f_{1,\theta}\|f_0)},
$$

where $D(f_{1,\theta}\|f_0)$ is the KL-divergence and $J_0$ the second-order moment of the log-likelihood between the (true) post-change and pre-change distribution, 
$$
D(f_{1,\theta}\|f_0) =\E_{1,\theta}\left[\log\frac{f_1(x_1,\theta)}{f_0(x_1)}\right],\quad J_0 = \E_{1,\theta}\left[\log^2\frac{f_1(x_1,\theta)}{f_0(x_1)}\right],
$$
$\hat D(f_{1,\theta}\|f_0) ,\hat J_0$ is defined similarly with true parameter $\theta$ replaced by the estimator,
$$\hat D(f_{1,\theta}\|f_0) = \E_{1,\theta}\left[\log\frac{f_1(x_t,\hat\theta_t)}{f_0(x_t)}\right] = D(f_{1,\theta}\|f_0) + O(1/w),$$
$$
\hat J_0 = \E_{1,\theta}\left[\log^2\frac{f_1(x_t,\hat\theta_t)}{f_0(x_t)}\right] = J_0 + O(1/w), \quad \text{for $t\geq w+1$}.
$$
\label{thm:edd_wl_cusum}
\end{theorem}
The computation complexity of this window-limited CUSUM depends on the parametric family and the optimization method to obtain the estimator $\hat \theta_t$. For the exponential family, $\hat\theta_t$ can be updated recursively in $O(1)$ time from $\hat\theta_{t-1}$ and for the general parametric family, this depends on the precision we want in solving optimization and the computation complexity may scale with window length $w$. Apart from this, the statistic $T_{t+1}^{\rm WL-CUSUM}$ can be updated in $O(1)$ time using the usual CUSUM recursion. Either way, it has better computation complexity than window-limited GLR because the latter one still needs to search for all possible change points within the much bigger sliding window.



\section{E-detectors: Nonparametric change-point detection}
\label{sec:non_parametric}

Recently, \cite{shin2022detectors} discuss a general framework based on the CUSUM procedure, which can be applied to non-parametric and composite pre- and post-change distribution families $\mathcal F_0$ and $\mathcal F_1$, named \textit{E-detectors}. A $\mathcal F_0$-E-detector is defined as a non-negative process $(M_t)_{t\geq 0}$ such that
$$
\E_{f,0}[M_\tau] \leq \E_{f,0}[\tau], \forall f\in\mathcal F_0, \forall \tau\in\mathcal T,
$$
where $\E_{f,0}$ is the expectation when the pre-change distribution is $f$ with no change-point, and $\mathcal T$ is the set of finite stopping times. It is easily seen by the optional stopping theorem that if $(M_t - t)_{t\geq 0}$ is a super-martingale under all probability measure in $\mathcal F_0$, then $(M_t)_{t\geq 0}$ must be an $\mathcal F_0$-E-detector. Let the stopping time
$$
\tau^{\rm E}(b) =\min\{t: M_t\geq \exp(b)\}. 
$$
Then the ARL is lower bounded by
\begin{equation}
\E_{f_0,0}[\tau^{\rm E}(b)]\geq \E_{f_0,0}[M_{\tau^{\rm E}(b)}]\geq \exp(b).
\label{eq:e_detector_arl}
\end{equation}
The classic CUSUM and SR procedures all satisfy the definition of E-detectors, and \eqref{eq:e_detector_arl} provides a fast and easy way to lower bound the ARL, which is good enough for showing the asymptotic optimality of the procedure.

For the CUSUM procedure, the final E-detector used for detecting an unknown anomaly is a mixture of baseline E-detectors. Each baseline E-detector is constructed via \textit{E-processes} $(M_t^{(k)})_{t\geq 0}$, $k = 1,2,3,\dots,$ which are generalizations of the likelihood ratio. Each E-process is responsible for detecting a change point at $k$. An example of the most basic setup described in Section \ref{sec:setup} is
\begin{equation}
M_t^{(k)} = \begin{cases}
\prod_{i=k}^t\frac{f_1(x_i)}{f_0(x_i)}& t\geq k,\\
1& 1\leq t<k.
\end{cases}
\label{eq:e-detector_example}
\end{equation}
In the non-parametric setup, the baseline E-process is in the form of
\begin{equation}
M_t^{(k)} = \begin{cases}
\prod_{i=k}^tL_i& t\geq k,\\
1& 1\leq t<k,
\end{cases}
\label{eq:e-process}
\end{equation}
where the multiplicative increments $L_i$ are non-negative and shared between all $k$ to keep the recursive rule \eqref{eq:simple_recursion}. $L_i$ also satisfies for any $f_0\in \mathcal F_0$, 
$
\E_{f_0}[L_i|x_1,\dots,x_{i-1}]\leq 1.
$
So the E-process $(M_t^{(k)})_{t\geq 0}$ is a non-negative super-martingale under the null hypothesis. An example of the multiplicative increments $L_i$ is by leveraging the concentration inequalities used to quantify the uncertainties in the (non-parametric) pre-change distribution $f_0$ (e.g. \cite{howard2020time}), which is in the form of 
$$
\E_{f_0}[\exp(\lambda s(x_i) - \varphi(\lambda)\nu(x_i))|x_1,\dots,x_{i-1}]\leq 1,\ \forall \lambda\in\Lambda
$$
for some known real function $s$, continuously differentiable convex function $\varphi$ and non-negative function $\nu$. Then the multiplicative increment is
\begin{equation}
L_i(\lambda) = \exp(\lambda s(x_i) - \varphi(\lambda)\nu(x_i)),
\label{eq:e_detector_increment}
\end{equation}
and the baseline E-detector is
\begin{equation}
M_t(\lambda) = \max_{1\leq k\leq t} M_t^{(k)}(\lambda), t = 0,1,2,\dots,
\label{eq:e-detector}
\end{equation}
where $\lambda$ remains to be chosen to optimize the EDD. 
To account for the unknown post-change distribution, a mixture of baseline E-detectors
$$
M_t = \sum_{i=1}^K \omega_iM_t(\lambda_i), \quad \sum_{i=1}^K \omega_i = 1
$$
can be taken so that the EDD for the worst possible anomaly will be small. Since the set of E-detectors is closed under convex operations, the mixture will still be an E-detector, and the ARL is controlled by \eqref{eq:e_detector_arl}. Let $\E_{1,f_1},\mathbb V_{1,f_1}$ represent the expectation and variance when the change starts with the first sample and the post-change distribution is $f_1$. Let $(L_1^{(\lambda)})_{\lambda\in\Lambda}$ be the set of all possible baseline increments in an E-process, and $\lambda^*$ be the parameter maximizing $\E_{1,f_1}[\log L_1^{(\lambda)}]$. Let $D(f_1\|\mathcal F_0) = \E_{1,f_1}[\log L_1^{(\lambda^*)}]$ be the divergence between $f_1$ and $\mathcal F_0$. Let $\Delta_U,\Delta_L$ be the upper and lower bound of the gap $\mathbb E_{1,f_1}[s(x_1)]/\mathbb E_{1,f_1}[\nu(x_1)]$ between $\mathcal F_0$ and all possible $f_1$, and let $\varphi^*$ be the convex conjugate of $\varphi.$ For i.i.d. data and a properly chosen mixture, the EDD is upper bounded by the following.
\begin{theorem}[EDD of a proper mixture](\cite[Theorem 4.3]{shin2022detectors})
For every $f_1\in\mathcal F_1$,
$$
\EDD(\tau^{\rm E}(b))\leq \frac{g_b}{D(f_1\|\mathcal F_0)} + \frac{\mathbb V_{1,f_1}(\log L_1^{(\lambda^*)})}{D(f_1\|\mathcal F_0)^2} + 1,
$$
where 
$$
g_b = \inf_{\eta>1} \eta\left(b + \log\left(1+ \left\lceil\log_\eta\frac{\varphi^*(\Delta_U)}{\varphi^*(\Delta_L)}\right\rceil\right)\right),
$$
given that there are sufficiently many baseline E-detectors in the mixture.
\end{theorem}
The first inequality is based on results in \cite{lorden1970excess}. Note that $D(f_1\|\mathcal F_0)$ is no longer the KL divergence, as there are no probability density functions. It is instead related to the construction of the set of baseline increments, in this case, the functions $s,v,\varphi$ in \eqref{eq:e_detector_increment} and $\Lambda$ the set of possible $\lambda$. So are the gap bounds $\Delta_U,\Delta_L$. Here $D(f_1\|\mathcal F_0)$ will be smaller than the KL divergence between the true $f_0$ and the true $f_1$, thus this procedure is not asymptotically optimal in the sense of Theorem \ref{thm:theoretical_lower_bound}. We believe no non-parametric methods can match this theoretical lower bound on EDD, as a cost of admitting better model robustness.

The computation and memory complexity is constant for each baseline E-detector, the same as the classic CUSUM procedure. The number of baseline E-detectors depends on the threshold $b$ for ARL control and the size of the gap between $\mathcal F_0$ and $\mathcal F_1$. For details, we refer to \cite{shin2022detectors}.

\section{Distribution-free Change-point detection procedures}
\label{sec:distribution_free}

Like has mentioned above, to our knowledge, existing distribution-free change-point detection procedures largely follow a similar strategy based on the Shewhart chart: scanning through different blocks of the data and finding the statistics $\{T_{t}\}_{t\geq 0}$ where a higher value suggests the distributions of the two blocks are different. Then the procedure declares a change-point when $T_t>b$ for some $t$. Such methods replace the log-likelihood ratio used in a parametric setup with alternative divergence, for example, the maximum mean discrepancy (MMD) \citep{li2019scan,song2022new}, kernel-based Fisher discriminant ratio \citep{harchaoui2008kernel}, similarity graphs \citep{chen2015graph, chu2019asymptotic}, marginal ranking \citep{lung2015homogeneity}, interpoint distance \citep{matteson2014nonparametric} to remove the strong assumption of known distribution family. Typically theoretical analysis for such procedures can be more challenging because of the lack of a handle through the likelihood ratio. Here we review a technique based on change-of-measure that can approximate the probability of extremely rare events leading to the ARL. We start from the analysis under the parametric setup (Section \ref{sec:extreme_in_random_fields}) and discuss its generalization to distribution-free change-point detection in Section \ref{sec:sub_distribution_free}. The complete and rigorous version is given in \cite{yakir2013extremes}. The EDD analysis remains an open question.

\subsection{ARL analysis via change-of-measure for extreme of random fields}
\label{sec:extreme_in_random_fields}
In this section, we review the ARL analysis of the GLR in \cite{siegmund1995using}. In this setup, the pre-change distribution $f_0\sim\mathcal N(0,1)$ and the post-change distribution $f_1\sim\mathcal N(\theta,1)$ where $\theta$ is unknown. The (window-limited) GLR statistic is in the form of 
$$
T_t^{\rm GLR} = \max_{min\{t-w,1\}\leq k\leq t}  \frac{\left(\sum_{i=k}^t x_i\right)^2}{2(t-k+1)},
$$
where for GLR, we use $w=  \infty$; note that the maximum over $\theta$ part is absorbed into the equation using a plug-in maximum likelihood estimator. 

For each $(t,k)$, we denote by $Z_{t,k} = \left(\sum_{i=k}^t x_i\right)/\sqrt{t-k+1}$ and each follows $N(0,1)$ under the null distribution. This kind of analysis can be generalized to setups where the log-likelihood ratios do not follow Gaussian distribution and have a possibly continuous-valued index $\theta$. Recall that a general change-point detection procedure is defined as the first time the detecting statistic exceeds threshold $b$, and in this case, we consider the one-sided stopping rule
\[
\tau^{\rm GLR} = \inf\left\{t: \max_{\min\{t-w,1\}\leq k\leq t} Z_{t,k} \geq \sqrt{2b}\right\}.
\]

\vspace{.1in}
\noindent
\textbf{Poisson approximation.} By examining the correlation between the $Z_{t,k}$s, it can be proved that $\tau^{\rm GLR}$ asymptotically follows exponential distribution \citep{siegmund1995using,yakir2013extremes}. Next we want to find out the mean by calculating $\mathbb P_{0} (\tau^{\rm GLR} \leq m)$ for $\log b<< m << b^{-1/2}e^b,$ 
\begin{equation}
\mathbb P_{0} (\tau^{\rm GLR} \leq m) = \mathbb P_{0} \left(\max_{\min\{t-w,1\}\leq k\leq  t\leq m} Z_{t,k} \geq \sqrt{2b}\right).
\label{eq:extreme_prob}
\end{equation}
And there would be $\ARL(\tau^{GLR}) \sim m/\mathbb P_0(\tau^{\rm GLR}\leq m)$, where the symbol $\sim$ means the ratio between the two sides converges to 1. 

\vspace{.1in}
\noindent\textbf{Decomposition on the local field of each $(t,k)$.} Note that the random field $\{Z_{t,k}\}$ is highly correlated because it is computed using highly overlapping data segments. For convenience, we denote the set of all possible $(t,k)$ by $S_{m,w} = \{(t,k)\in\mathbb N^2: \min\{t-w,1\}\leq k\leq t\leq m\}$ and the event $\{\max_{(t,k)\in S_{m,w}}  Z_{t,k} \geq \sqrt{2b}\}$ by $ A_{m,w}$. To analyze the tail probability on the right-hand-side of \eqref{eq:extreme_prob}, we decompose it as the sum of the probability $Z_{t,k}\geq \sqrt{2b}$ for each $(t,k)\in  S_{m,k}$ while also take into account their correlation. Let $\boldsymbol 1$ be the indicator function. For any function $g:\mathbb R\to\mathbb R$,
\begin{align*}&\ \mathbb P_{0}  \left(\max_{(t,k)\in  S_{m,w}} Z_{t,k} \geq \sqrt{2b}\right)\\
=&\ \E_0\left[\frac{\sum_{(t,k)\in  S_{m,w}}\exp(g(Z_{t,k}))}{\sum_{(t',k')\in  S_{m,w}}\exp(g(Z_{t',k'}))}\boldsymbol 1\{  A_{m,w}\}\right]\\
=&\ \sum_{(t,k)\in  S_{m,w}}\E_{0}\left[\frac{\exp(g(Z_{t,k}))}{\sum_{(t',k')\in  S_{m,w}}\exp(g(Z_{t',k'}))}\boldsymbol 1\{  A_{m,w}\}\right]
\end{align*}

\noindent\textbf{Exponential tilting and change-of-measure.}
We use a change-of-measure technique similar to that has been used in  Section \ref{sec:performance_analysis} equation \eqref{eq:change_of_measure} so that under the new probability measure, $\{Z_{t,k}\geq \sqrt{2b}\}$ happens with a much higher rate. In this setup for each $(t,k)$ we change to the measure $\mathbb P_{t,k}$ and its expectation $\E_{t,k}$, under which $x_k,x_{k+1},\dots,x_t \sim \mathcal N(\sqrt{2b}/\sqrt{t-k+1},1)$ while the rest of the samples still follow $\mathcal N(0,1)$. So now $Z_{t,k}\sim \mathcal N(\sqrt{2b},1)$ under $\mathbb P_{t,k}$, and the event $  A_{m,w}$ happens with high probability. Let $g(Z_{t,k}) = \sqrt{2b}Z_{t,k} - b$ be the log-likelihood ratio between $\mathbb P_{t,k}$ and $\mathbb P_0$. When $Z_{t,k}$ is not normally distributed, exponential tilting is done by letting $g(z) = \gamma z - \varphi(\gamma)$, where $\varphi$ is the log moment generating function of $Z_{t,k}$ under the null distribution. Then $\E_0[\exp(g(Z_{t,k}))] = 1$ and it can be regarded as the likelihood ratio between $\mathbb P_0$ and some distribution $\mathbb P_{t,k}$. By choosing a proper $\gamma$ such that $\E_{t,k}[Z_{t,k}]$ equals to the threshold $b$, $A_{m,w}$ happens with high probability.
\begin{align*}
    &\ \mathbb P_{0}  \left(\max_{(t,k)\in  S_{m,w}} Z_{t,k} \geq \sqrt{2b}\right)\\
    =&\ \sum_{(t,k)\in  S_{m,w}}\E_{t,k}\left[\frac{1}{\sum_{(t',k')\in  S_{m,w}}\exp(\sqrt{2b}Z_{t',k'}-b)}\boldsymbol 1\{  A_{m,w}\}\right].
\end{align*}
Next step, we decompose the expectation into a local and a global term which are ``asymptotically independent,'' the meaning of which will be explained later.
\begin{align*}
&\ \sum_{(t,k)\in  S_{m,w}} \E_{t,k}\Bigg[\frac{1}{\sum_{(t',k')\in  S_{m,w}}\exp(\sqrt{2b}Z_{t',k'}-b)}\boldsymbol 1\{  A_{m,w}\}\Bigg]\\
=&\ \sum_{(t,k)\in  S_{m,w}}\E_{t,k}\Bigg[\frac{\max_{(t',k')\in  S_{m,w}}\exp(\sqrt{2b}(Z_{t',k'} - Z_{t,k}))}{\sum_{(t',k')\in  S_{m,w}}\exp(\sqrt{2b}(Z_{t',k'} - Z_{t,k}))}\exp\left(-\max_{(t',k')\in  S_{m,w}}(\sqrt{2b}Z_{t',k'}-b)\right)\boldsymbol 1\{ A_{m,w}\}\Bigg].
\end{align*}
For convenience, let $$
S_{t,k} = \sum_{(t',k')\in  S_{m,w}}\exp(\sqrt{2b}(Z_{t',k'} - Z_{t,k})),\ 
M_{t,k} = \max_{(t',k')\in  S_{m,w}}\exp(\sqrt{2b}(Z_{t',k'} - Z_{t,k})).
$$
Now after replacing certain terms with $S_{t,k}$ and $M_{t,k}$,
\begin{equation}
\mathbb P(  A_{m,w}) = \sum_{(t,k)\in  S_{m,w}}\E_{t,k}\left[\frac{M_{t,k}}{S_{t,k}}\exp\left(-(\log M_{t,k} +\sqrt{2b} Z_{t,k} - b)\right)\boldsymbol 1\{\log M_{t,k} + \sqrt{2b}Z_{t,k} -b \geq b\}\right].
\label{eq:prob_m}
\end{equation}
Next we discuss the expectation as 
\begin{equation}
b\to\infty,t-k\to\infty\text{ and }b/(t-k+1) \to c \text{ for some constant $c$.}
\label{eq:asymptotic_regime}
\end{equation}

\noindent\textbf{The Local term.} The first term $M_{t,k}/S_{t,k}$ in \eqref{eq:prob_m} for every $(t,k)\in S_{m,w}$, or the tuple $(M_{t,k},S_{t,k})$, is a local random variable meaning it can be well approximated with high probability using values in the random field $(Z_{t',k'})_{(t',k')\in S_{m,w}}$ which index is constantly close to $(t,k)$. And $(M_{t,k},S_{t,k})$ converge in distribution to some local random variables which we'll define later. Since it is a Gaussian random field with mean 0 under $\mathbb P_0$, the local properties are decided by the covariance under the asymptotic regime \eqref{eq:asymptotic_regime}. For fixed $i,j,i',j'\in\mathbb Z$, it can be verified that
\begin{equation}
\Cov_0[Z_{t+i,k-j},Z_{t+i',k-j'}] = 1 - \frac{c}{2b}|i-i'|-\frac{c}{2b}|j-j'| + o(b^{-1}).
\label{eq:local_covariance}
\end{equation}
And consequently,
$$
\E_{t,k}[\sqrt{2b}(Z_{t+i,k-j} - Z_{t,k})] = -c|i|-c|j| + o(1),
$$
\begin{align*}
\Cov_{t,k}[\sqrt{2b}(Z_{t+i,k-j}-Z_{t,k}),&\ \sqrt{2b}(Z_{t+i',k-j'} - Z_{t,k})] \\
= &\ 2c\boldsymbol{1}\{ii'>0\}\min\{|i|,|i'|\} + 2c \boldsymbol{1}\{jj'>0\}\min\{|j|,|j'|\}.
\end{align*}
The local random field $\{\sqrt{2b}(Z_{t+i,k-j}-Z_{t,k})\}_{i,j}$ under $\mathbb P_{t,k}$ converge in distribution to \begin{equation}
\{Y_i + Y'_j\}_{i,j},
\label{eq:two_gaussian_random_walk}
\end{equation}
where $\{Y_i\}_{i\in\mathbb Z}$ is a two-sided Gaussian random walk with a negative drift $c$ and variance $2c$, i.e.
$$\E[Y_i] = -c|i|, \quad \Cov[Y_i,Y_{i'}] = 2c\boldsymbol 1\{ii'>0\}\min\{|i|,|i'|\}.
$$
$\{Y'_j\}_{j\in\mathbb Z}$ is an i.i.d. copy of $\{Y_i\}_{i\in\mathbb Z}$. Let 
$$
\mathcal M_c = \max_{i\in\mathbb Z} \exp(Y_i),\quad \mathcal S_c = \sum_{i\in\mathbb Z} \exp(Y_i),
$$
The expectation $\E[\mathcal M_c/\mathcal S_c]$ is called the Mill's ratio (see Chapter 2.2 in \cite{yakir2013extremes}) and it has the following form.
\begin{equation}
\E[\mathcal M_c/\mathcal S_c] = c\nu(\sqrt{2c}),
\label{eq:mills_ratio}
\end{equation}
where the special function
$$
\nu(x) = 2x^{-2}\exp\left(-2\sum_{i=1}^\infty i^{-1}\Phi(-xi^{1/2}/2)\right),
$$
is closely related to the Laplace transform of the overshoot over the boundary of a random walk. A numerical approximation to the Mill's ratio is
$$
\nu(x) \approx \frac{(2/x)(\Phi(x/2)-0.5)}{(x/2)\Phi(x/2) + \phi(x/2)},
$$
here $\phi,\Phi$ is the p.d.f. and c.d.f. of the standard normal distribution.

Back to the GLR analysis, it can be shown that $(M_{t,k},S_{t,k})$ converge in distribution to $(\mathcal M_c\mathcal M'_c,\mathcal S_c\mathcal S'_c)$ where $(\mathcal M'_c,\mathcal S'_c)$ is an i.i.d. copy of $(\mathcal M_c,\mathcal S_c)$. This is because our local random field is two-dimensional and asymptotically we can break it down using \eqref{eq:two_gaussian_random_walk} into the sum of two independent Gaussian random walks. Using the Mill's ratio, the expectation $\E_{t,k}[M_{t,k}/S_{t,k}]$ over the two-dimensional random field is
$$
\E_{t,k}[ M_{t,k}/ S_{t,k}] \to \E[\mathcal M_c/\mathcal S_c]^2 = c^2\nu^2(\sqrt{2c}).
$$

\noindent\textbf{The global term.} The second term $\exp\left(-(\log M_{t,k} + \sqrt{2b}Z_{t,k}-b)\right)\boldsymbol 1\{\log M_{t,k} + \sqrt{2b}Z_{t,k}-b \geq b\}$ in \eqref{eq:prob_m} for every $(t,k)\in S_{m,w}$ is considered a global one because asymptotically it is not relevant with the constantly many observations around time $t$ and $k$, but only those in between. Under the asymptotic regime \eqref{eq:asymptotic_regime}, $\sqrt{2b}Z_{t,k} - b\sim \mathcal N(b,2b)$ but the distribution of $M_{t,k}$ converge to $\mathcal M_c\mathcal M'_c$ which does not scale with $b$. It can be shown that the distribution of $\log M_{t,k} + \sqrt{2b}Z_{t,k}-b$ together is dominated by $\sqrt{2b}Z_{t,k}-b$ and it is so close to $\mathcal N(b,2b)$ that for each $(t,k)$ in \eqref{eq:prob_m} the global term can be treated as such and is `independent' of the local term. And the expectation of the global term converges to
$$
\E_{t,k}\left[\exp\left(-(\sqrt{2b}Z_{t,k}-b)\right)\boldsymbol 1\{\sqrt{2b}Z_{t,k}-b \geq b\}\right] \to \int_b^\infty \frac{1}{\sqrt{4\pi b}}\exp(-x)dx = e^{-b}/\sqrt{4\pi b}
$$
because the probability density function of $\sqrt{2b}Z_{t,k}-b$ is almost always $1/\sqrt{4\pi b}$ around $b$.

\noindent Together, we know that asymptotically the expectation in \eqref{eq:prob_m} is
\begin{equation}
\begin{split}
\E_{t,k}\bigg[\frac{M_{t,k}}{S_{t,k}}\exp\left(-(\log M_{t,k} +\sqrt{2b} Z_{t,k} - b\right) \boldsymbol 1\{\log M_{t,k} + &\ \sqrt{2b}Z_{t,k} - b \geq b\}\bigg]\\
&\sim \frac{e^{-b}}{\sqrt{4\pi b}}c^2\nu^2(\sqrt{2c}).
\end{split}
\label{eq:decomposed_prob}
\end{equation}

\noindent \textbf{The ARL.} After checking several good conditions on $\E_{t,k}[\cdot]$ (e.g., uniform integrability), we can now calculate the probability $\mathbb P(A_{m,w})$ by substituting the summation with integral. Assume as $b\to\infty$, the window length $w$ satisfies $b/w\to C$, and $m>>w$,
\begin{align*}
\mathbb P(A_{m,w}) \sim &\ \frac{e^{-b}}{\sqrt{4\pi b}}\sum_{(t,k)\in S_{m,w}}\E[\mathcal M_{b/(t-k+1)}/\mathcal S_{b/t-k+1}]^2\\
\sim &\ \frac{e^{-b}}{\sqrt{4\pi b}} b\int_0^{C^{-1}} m\E[\mathcal M_{x^{-1}}/\mathcal S_{x^{-1}}]^2dx,
\end{align*}
where we take $x = (t-k+1)/b$ when changing to integral. By Poisson approximation,
$$
\ARL(\tau^{\rm GLR}) \sim \sqrt{4\pi}b^{-0.5} e^{b}\left(\int_0^{C^{-1}} x^{-2}\nu^2(\sqrt{2/x})dx\right)^{-1}.
$$

\subsection{Scan statistic via kernel-based MMD}
\label{sec:sub_distribution_free}

In \citep{li2019scan}, the authors propose to use the unbiased kernel-based maximum mean divergence (MMD) statistic to measure the difference between two blocks of samples. For $X = \{x_1,x_2,\dots,x_n\}$ and $Y = \{y_1,y_2,\dots,y_n\}$, the unbiased MMD is
$$
\mathrm{MMD}(X,Y) = \frac{1}{n(n-1)}\sum_{i\neq j}(k(x_i,x_j) + k(y_i,y_j) - k(x_i,y_j) - k(y_i,x_j)),
$$
where $k$ is a kernel in the reproducing kernel Hilbert space. The expectation of such measure is 0 only when samples in $X$ and $Y$ share the same distribution; otherwise, the expectation is strictly positive. Assuming the abundance of reference data, at every time $t$ the average unbiased MMD between the most recent $B$ samples and $N$ reference blocks is computed, and then normalized to get the test statistic $T_t^{\rm MMD}$:
$$
Z_t = \frac{1}{N}\sum_{i=1}^N \mathrm{MMD_u}(X_{t},Y_{t,i}), \quad\quad  T_t^{\rm MMD} = \frac{Z_t}{\mathbb V_0(Z_t)^{1/2}}.
$$
Here $X_t$ is the collection of the most recent $B$ samples, and the reference blocks $Y_{t,i},i=1,\dots,N$ is formed by taking $NB$ samples without replacement from reference data. The variance $\mathbb V_0(Z_t)^{1/2}$ can be pre-computed from reference data. And procedure raises an alarm at $\tau^{\rm MMD} = \min\{t:T_t^{\rm MMD}>b\}$. 

During each update, the statistic $T_{t+1}^{\rm MMD}$ can be computed recursively from $T_{t}^{\rm MMD}$. This is because the reference block $Y_{t+1,i}$ is formed by adding a new sample from reference data to $Y_{t,i}$ and removing the oldest one, for $i=1,\dots,N$. $X_{t+1}$ is formed using the same way. For convenience, let $$
h(x,x',y,y') = k(x,x') + k(y,y') - k(x,y') + k(y,x'),
$$
and let $Y_{t,i} = (y_{t-B+1,i},\dots,y_{t,i}). $ for every $i=1,\dots,N$. Then $T_{t+1}^{\rm MMD}$ can be updated from $T_t^{\rm MMD}$ by adding and subtracting the values $h(x_s,x_{t+1},y_{s,i},y_{t+1,i}),$ $h(x_s,x_{t-B+2},y_{s,i},y_{t-B+2,i})$ $s=t-B+2,\dots,t, i=1,\dots,N$. This leads to a computation and memory complexity of $O(NB)$.

The ARL analysis of such statistics follows the example given in Section \ref{sec:extreme_in_random_fields}. The covariance between the statistics $\{T_t\}_{t\geq 0}$ is provided in a close form and can be computed easily from reference data. With a proper asymptotic regime, the local covariance can be found in the form of \eqref{eq:local_covariance}. The rest of the analysis carries through as if $\{T_t\}_{t\geq 0}$ follows Gaussian distribution.

Of course, the above analysis can be less accurate when $\{Z_t\}_{t\geq 0}$ does not converge to a Gaussian distribution. Skewness correction can be performed to account for the non-zero third-order moment. This method first appears in \cite{tu1999maximum} and then is modified in \cite{tang2001mapping} and becomes useful in analyzing the ARL or false alarm rate in many change-point detection procedures based on scanning statistics.

\section{Conclusion}
\label{sec:conclusion}

In this paper, we review the classic procedures in sequential change-point detection featuring the computation versus performance trade-off and the classic performance analysis methods, which are still quite useful in analyzing new change-point detection procedures and other problems (for details, we refer to \cite{yakir2013extremes}). We also discuss several recent papers covering computation and model robustness considerations. For completeness, a longer survey on the extensions and modern applications of sequential change-point detection can be found in \cite{xie2021sequential}. Finally, another important direction we have not mentioned is machine learning algorithm-based sequential change-point detection to handle situations such as missing data \citep{xie2012change,londschien2021change}, the existence of outliers \citep{fearnhead2019changepoint}, and cost-sensitive active sampling \citep{gundersen2021active}. 

\section*{Acknowledgements}

This work is partially supported by an NSF CAREER CCF-1650913, NSF DMS-2134037, CMMI-2015787, CMMI-2112533, DMS-1938106, and DMS-1830210, and a Coca-Cola Foundation fund.

\bibliographystyle{apalike}
\bibliography{references}

\begin{thebibliography}{}

\bibitem[Abbasi and Haq, 2019]{abbasi2019optimal}
Abbasi, S. and Haq, A. (2019).
\newblock Optimal {CUSUM} and adaptive {CUSUM} charts with auxiliary
  information for process mean.
\newblock {\em Journal of Statistical Computation and Simulation},
  89(2):337--361.

\bibitem[Basseville et~al., 1993]{basseville1993detection}
Basseville, M., Nikiforov, I.~V., et~al. (1993).
\newblock {\em Detection of abrupt changes: theory and application}, volume
  104.
\newblock Prentice-Hall Englewood Cliffs.

\bibitem[Cao et~al., 2018]{cao2018sequential}
Cao, Y., Xie, L., Xie, Y., and Xu, H. (2018).
\newblock Sequential change-point detection via online convex optimization.
\newblock {\em Entropy}, 20(2):108.

\bibitem[Casella and Berger, 2021]{casella2021statistical}
Casella, G. and Berger, R.~L. (2021).
\newblock {\em Statistical inference}.
\newblock Cengage Learning.

\bibitem[Chen and Zhang, 2015]{chen2015graph}
Chen, H. and Zhang, N. (2015).
\newblock Graph-based change-point detection.
\newblock {\em The Annals of Statistics}, 43(1):139--176.

\bibitem[Chen et~al., 2020]{chen2017textsf}
Chen, J., Kim, S.-H., and Xie, Y. (2020).
\newblock $\textsf{S}^3t$: An efficient score-statistic for spatio-temporal
  surveillance.
\newblock {\em Sequential Analysis}, 39:563--592.

\bibitem[Chen et~al., 2022]{chen2022high}
Chen, Y., Wang, T., and Samworth, R.~J. (2022).
\newblock High-dimensional, multiscale online changepoint detection.
\newblock {\em Journal of the Royal Statistical Society Series B: Statistical
  Methodology}, 84(1):234--266.

\bibitem[Chen et~al., 2015]{chen2015quickest}
Chen, Y.~C., Banerjee, T., Dominguez-Garcia, A.~D., and Veeravalli, V.~V.
  (2015).
\newblock Quickest line outage detection and identification.
\newblock {\em IEEE Transactions on Power Systems}, 31(1):749--758.

\bibitem[Chu and Chen, 2019]{chu2019asymptotic}
Chu, L. and Chen, H. (2019).
\newblock Asymptotic distribution-free change-point detection for multivariate
  and non-euclidean data.
\newblock {\em The Annals of Statistics}, 47(1):382--414.

\bibitem[Dehning et~al., 2020]{dehning2020inferring}
Dehning, J., Zierenberg, J., Spitzner, F.~P., Wibral, M., Neto, J.~P., Wilczek,
  M., and Priesemann, V. (2020).
\newblock Inferring change points in the spread of covid-19 reveals the
  effectiveness of interventions.
\newblock {\em Science}, 369(6500):eabb9789.

\bibitem[Fearnhead and Rigaill, 2019]{fearnhead2019changepoint}
Fearnhead, P. and Rigaill, G. (2019).
\newblock Changepoint detection in the presence of outliers.
\newblock {\em Journal of the American Statistical Association},
  114(525):169--183.

\bibitem[Gundersen et~al., 2021]{gundersen2021active}
Gundersen, G.~W., Cai, D., Zhou, C., Engelhardt, B.~E., and Adams, R.~P.
  (2021).
\newblock Active multi-fidelity {Bayesian} online changepoint detection.
\newblock In {\em Uncertainty in Artificial Intelligence}, pages 1916--1926.
  PMLR.

\bibitem[Hallac et~al., 2015]{hallac2015network}
Hallac, D., Leskovec, J., and Boyd, S. (2015).
\newblock Network lasso: Clustering and optimization in large graphs.
\newblock In {\em Proceedings of the 21th ACM SIGKDD international conference
  on knowledge discovery and data mining}, pages 387--396.

\bibitem[Harchaoui et~al., 2008]{harchaoui2008kernel}
Harchaoui, Z., Moulines, E., and Bach, F. (2008).
\newblock Kernel change-point analysis.
\newblock {\em Advances in neural information processing systems}, 21.

\bibitem[Howard et~al., 2020]{howard2020time}
Howard, S.~R., Ramdas, A., McAuliffe, J., and Sekhon, J. (2020).
\newblock Time-uniform {Chernoff} bounds via nonnegative supermartingales.
\newblock {\em Probability Surveys}, 17:257--317.

\bibitem[Krieger et~al., 1999]{krieger1999detecting}
Krieger, A.~M., Pollak, M., and Yakir, B. (1999).
\newblock Detecting a change in regression: first-order optimality.
\newblock {\em The Annals of Statistics}, 27(6):1896--1913.

\bibitem[Lai, 1998]{lai1998information}
Lai, T.~L. (1998).
\newblock Information bounds and quick detection of parameter changes in
  stochastic systems.
\newblock {\em IEEE Transactions on Information theory}, 44(7):2917--2929.

\bibitem[Lai and Shan, 1999]{lai1999efficient}
Lai, T.~L. and Shan, J.~Z. (1999).
\newblock Efficient recursive algorithms for detection of abrupt changes in
  signals and control systems.
\newblock {\em IEEE Transactions on Automatic Control}, 44(5):952--966.

\bibitem[Lakhina et~al., 2004]{lakhina2004diagnosing}
Lakhina, A., Crovella, M., and Diot, C. (2004).
\newblock Diagnosing network-wide traffic anomalies.
\newblock {\em ACM SIGCOMM computer communication review}, 34(4):219--230.

\bibitem[Lee and Kriegman, 2005]{lee2005online}
Lee, K.-C. and Kriegman, D. (2005).
\newblock Online learning of probabilistic appearance manifolds for video-based
  recognition and tracking.
\newblock In {\em 2005 IEEE Computer Society Conference on Computer Vision and
  Pattern Recognition (CVPR'05)}, volume~1, pages 852--859. IEEE.

\bibitem[Li et~al., 2019]{li2019scan}
Li, S., Xie, Y., Dai, H., and Song, L. (2019).
\newblock Scan {B}-statistic for kernel change-point detection.
\newblock {\em Sequential Analysis}, 38(4):503--544.

\bibitem[Londschien et~al., 2021]{londschien2021change}
Londschien, M., Kov{\'a}cs, S., and B{\"u}hlmann, P. (2021).
\newblock Change-point detection for graphical models in the presence of
  missing values.
\newblock {\em Journal of Computational and Graphical Statistics},
  30(3):768--779.

\bibitem[Lorden, 1970]{lorden1970excess}
Lorden, G. (1970).
\newblock On excess over the boundary.
\newblock {\em The Annals of Mathematical Statistics}, 41(2):520--527.

\bibitem[Lorden, 1971]{lorden1971procedures}
Lorden, G. (1971).
\newblock Procedures for reacting to a change in distribution.
\newblock {\em The Annals of Mathematical Statistics}, pages 1897--1908.

\bibitem[Lorden and Eisenberger, 1973]{lorden1973detection}
Lorden, G. and Eisenberger, I. (1973).
\newblock Detection of failure rate increases.
\newblock {\em Technometrics}, 15(1):167--175.

\bibitem[Lorden and Pollak, 2005]{lorden2005nonanticipating}
Lorden, G. and Pollak, M. (2005).
\newblock Non-anticipating estimation applied to sequential analysis and
  changepoint detection.
\newblock {\em The Annals of Statistics}, 33(3):1422--1454.

\bibitem[Lucas, 1982]{lucas1982combined}
Lucas, J.~M. (1982).
\newblock Combined {Shewhart-CUSUM} quality control schemes.
\newblock {\em Journal of Quality Technology}, 14(2):51--59.

\bibitem[Lung-Yut-Fong et~al., 2015]{lung2015homogeneity}
Lung-Yut-Fong, A., L{\'e}vy-Leduc, C., and Capp{\'e}, O. (2015).
\newblock Homogeneity and change-point detection tests for multivariate data
  using rank statistics.
\newblock {\em Journal de la Soci{\'e}t{\'e} Fran{\c{c}}aise de Statistique},
  156(4):133--162.

\bibitem[Magesh et~al., 2022]{magesh2022multiple}
Magesh, A., Veeravalli, V.~V., Roy, A., and Jha, S. (2022).
\newblock Multiple testing framework for out-of-distribution detection.
\newblock {\em arXiv preprint arXiv:2206.09522}.

\bibitem[Maillard, 2019]{maillard2019sequential}
Maillard, O.-A. (2019).
\newblock Sequential change-point detection: Laplace concentration of scan
  statistics and non-asymptotic delay bounds.
\newblock In Garivier, A. and Kale, S., editors, {\em Proceedings of the 30th
  International Conference on Algorithmic Learning Theory}, volume~98 of {\em
  Proceedings of Machine Learning Research}, pages 610--632. PMLR.

\bibitem[Malladi et~al., 2013]{malladi2013online}
Malladi, R., Kalamangalam, G.~P., and Aazhang, B. (2013).
\newblock Online {Bayesian} change point detection algorithms for segmentation
  of epileptic activity.
\newblock In {\em 2013 Asilomar conference on signals, systems and computers},
  pages 1833--1837. IEEE.

\bibitem[Matteson and James, 2014]{matteson2014nonparametric}
Matteson, D.~S. and James, N.~A. (2014).
\newblock A nonparametric approach for multiple change point analysis of
  multivariate data.
\newblock {\em Journal of the American Statistical Association},
  109(505):334--345.

\bibitem[Mei, 2006]{mei2006sequential}
Mei, Y. (2006).
\newblock Sequential change-point detection when unknown parameters are present
  in the pre-change distribution.
\newblock {\em The Annals of Statistics}, 34(1):92--122.

\bibitem[Moustakides, 1986]{moustakides1986optimal}
Moustakides, G.~V. (1986).
\newblock Optimal stopping times for detecting changes in distributions.
\newblock {\em the Annals of Statistics}, 14(4):1379--1387.

\bibitem[Moustakides, 2014]{moustakides2014multiple}
Moustakides, G.~V. (2014).
\newblock Multiple optimality properties of the {Shewhart} test.
\newblock {\em Sequential Analysis}, 33(3):318--344.

\bibitem[Page, 1954]{page1954continuous}
Page, E.~S. (1954).
\newblock Continuous inspection schemes.
\newblock {\em Biometrika}, 41(1/2):100--115.

\bibitem[Pavlov, 1991]{pavlov1991sequential}
Pavlov, I.~V. (1991).
\newblock Sequential procedure of testing composite hypotheses with
  applications to the {Kiefer–Weiss} problem.
\newblock {\em Theory of Probability \& Its Applications}, 35(2):280--292.

\bibitem[Peel and Clauset, 2015]{peel2015detecting}
Peel, L. and Clauset, A. (2015).
\newblock Detecting change points in the large-scale structure of evolving
  networks.
\newblock In {\em Twenty-Ninth AAAI Conference on Artificial Intelligence}.

\bibitem[Pollak, 1985]{pollak1985optimal}
Pollak, M. (1985).
\newblock Optimal detection of a change in distribution.
\newblock {\em The Annals of Statistics}, pages 206--227.

\bibitem[Pollak and Krieger, 2013]{Pollak2013}
Pollak, M. and Krieger, A.~M. (2013).
\newblock Shewhart revisited.
\newblock {\em Sequential Analysis}, 32(2):230--242.

\bibitem[Pollak and Siegmund, 1991]{pollak1991sequential}
Pollak, M. and Siegmund, D. (1991).
\newblock Sequential detection of a change in a normal mean when the initial
  value is unknown.
\newblock {\em The Annals of Statistics}, 19(1):394--416.

\bibitem[Pollak and Tartakovsky, 2009]{pollak2009optimality}
Pollak, M. and Tartakovsky, A.~G. (2009).
\newblock Optimality properties of the {Shiryaev-Roberts} procedure.
\newblock {\em Statistica Sinica}, pages 1729--1739.

\bibitem[Polunchenko and Tartakovsky, 2010]{polunchenko2010optimality}
Polunchenko, A.~S. and Tartakovsky, A.~G. (2010).
\newblock On optimality of the {Shiryaev-Roberts} procedure for detecting a
  change in distribution.
\newblock {\em The Annals of Statistics}, 38(6):3445--3457.

\bibitem[Polunchenko and Tartakovsky, 2012]{polunchenko2012state}
Polunchenko, A.~S. and Tartakovsky, A.~G. (2012).
\newblock State-of-the-art in sequential change-point detection.
\newblock {\em Methodology and computing in applied probability},
  14(3):649--684.

\bibitem[Raghavan and Veeravalli, 2010]{raghavan2010quickest}
Raghavan, V. and Veeravalli, V.~V. (2010).
\newblock Quickest change detection of a {Markov} process across a sensor
  array.
\newblock {\em IEEE Transactions on Information Theory}, 56(4):1961--1981.

\bibitem[Raginsky et~al., 2012]{raginsky2012sequential}
Raginsky, M., Willett, R.~M., Horn, C., Silva, J., and Marcia, R.~F. (2012).
\newblock Sequential anomaly detection in the presence of noise and limited
  feedback.
\newblock {\em IEEE Transactions on Information Theory}, 58(8):5544--5562.

\bibitem[Ren et~al., 2019]{ren2019likelihood}
Ren, J., Liu, P.~J., Fertig, E., Snoek, J., Poplin, R., Depristo, M., Dillon,
  J., and Lakshminarayanan, B. (2019).
\newblock Likelihood ratios for out-of-distribution detection.
\newblock {\em Advances in neural information processing systems}, 32.

\bibitem[Robbins and Siegmund, 1972]{robbins1972class}
Robbins, H. and Siegmund, D. (1972).
\newblock A class of stopping rules for testing parametric hypotheses.
\newblock In {\em Proc. Sixth Berkeley Symp. Math. Statist. Probab}, volume~4,
  pages 37--41.

\bibitem[Robbins and Siegmund, 1974]{robbins1974expected}
Robbins, H. and Siegmund, D. (1974).
\newblock The expected sample size of some tests of power one.
\newblock {\em The Annals of Statistics}, 2(3):415--436.

\bibitem[Roberts, 1966]{roberts1966comparison}
Roberts, S. (1966).
\newblock A comparison of some control chart procedures.
\newblock {\em Technometrics}, 8(3):411--430.

\bibitem[Romano et~al., 2023]{romano2023fast}
Romano, G., Eckley, I.~A., Fearnhead, P., and Rigaill, G. (2023).
\newblock Fast online changepoint detection via functional pruning {CUSUM}
  statistics.
\newblock {\em Journal of Machine Learning Research}, 24(81):1--36.

\bibitem[Shewhart, 1925]{shewhart1925application}
Shewhart, W.~A. (1925).
\newblock The application of statistics as an aid in maintaining quality of a
  manufactured product.
\newblock {\em Journal of the American Statistical Association},
  20(152):546--548.

\bibitem[Shewhart, 1931]{shewhart1931economic}
Shewhart, W.~A. (1931).
\newblock {\em Economic control of quality of manufactured product}.
\newblock Macmillan And Co Ltd, London.

\bibitem[Shin et~al., 2022]{shin2022detectors}
Shin, J., Ramdas, A., and Rinaldo, A. (2022).
\newblock E-detectors: a nonparametric framework for online changepoint
  detection.
\newblock {\em arXiv preprint arXiv:2203.03532}.

\bibitem[Shiryaev, 1963]{shiryaev1963optimum}
Shiryaev, A.~N. (1963).
\newblock On optimum methods in quickest detection problems.
\newblock {\em Theory of Probability \& Its Applications}, 8(1):22--46.

\bibitem[Siegmund, 1985]{siegmund1985sequential}
Siegmund, D. (1985).
\newblock {\em Sequential analysis: tests and confidence intervals}.
\newblock Springer Science \& Business Media.

\bibitem[Siegmund and Venkatraman, 1995]{siegmund1995using}
Siegmund, D. and Venkatraman, E. (1995).
\newblock Using the generalized likelihood ratio statistic for sequential
  detection of a change-point.
\newblock {\em The Annals of Statistics}, pages 255--271.

\bibitem[Song and Chen, 2022]{song2022new}
Song, H. and Chen, H. (2022).
\newblock New kernel-based change-point detection.
\newblock {\em arXiv preprint arXiv:2206.01853}.

\bibitem[Sparks, 2000]{sparks2000cusum}
Sparks, R.~S. (2000).
\newblock {CUSUM} charts for signalling varying location shifts.
\newblock {\em Journal of Quality Technology}, 32(2):157--171.

\bibitem[Tang and Siegmund, 2001]{tang2001mapping}
Tang, H.-K. and Siegmund, D. (2001).
\newblock Mapping quantitative trait loci in oligogenic models.
\newblock {\em Biostatistics}, 2(2):147--162.

\bibitem[Tartakovsky et~al., 2014]{tartakovsky2014sequential}
Tartakovsky, A., Nikiforov, I., and Basseville, M. (2014).
\newblock {\em Sequential analysis: Hypothesis testing and changepoint
  detection}.
\newblock CRC Press.

\bibitem[Tartakovsky, 2005]{tartakovsky2005asymptotic}
Tartakovsky, A.~G. (2005).
\newblock Asymptotic performance of a multichart {CUSUM} test under false alarm
  probability constraint.
\newblock In {\em Proceedings of the 44th IEEE Conference on Decision and
  Control}, pages 320--325. IEEE.

\bibitem[Tartakovsky, 2014]{tartakovsky2014nearly}
Tartakovsky, A.~G. (2014).
\newblock Nearly optimal sequential tests of composite hypotheses revisited.
\newblock {\em Proceedings of the Steklov Institute of Mathematics},
  287(1):268--288.

\bibitem[Tartakovsky, 2019]{tartakovsky2019sequential}
Tartakovsky, A.~G. (2019).
\newblock {\em Sequential change detection and hypothesis testing: general
  non-iid stochastic models and asymptotically optimal rules}.
\newblock Chapman and Hall/CRC.

\bibitem[Tartakovsky et~al., 2012a]{tartakovsky2012third}
Tartakovsky, A.~G., Pollak, M., and Polunchenko, A.~S. (2012a).
\newblock Third-order asymptotic optimality of the generalized
  {Shiryaev-Roberts} changepoint detection procedures.
\newblock {\em Theory of Probability \& Its Applications}, 56(3):457--484.

\bibitem[Tartakovsky et~al., 2012b]{tartakovsky2012efficient}
Tartakovsky, A.~G., Polunchenko, A.~S., and Sokolov, G. (2012b).
\newblock Efficient computer network anomaly detection by changepoint detection
  methods.
\newblock {\em IEEE Journal of Selected Topics in Signal Processing},
  7(1):4--11.

\bibitem[Truong et~al., 2020]{truong2020selective}
Truong, C., Oudre, L., and Vayatis, N. (2020).
\newblock Selective review of offline change point detection methods.
\newblock {\em Signal Processing}, 167:107299.

\bibitem[Tu and Siegmund, 1999]{tu1999maximum}
Tu, I.-P. and Siegmund, D. (1999).
\newblock The maximum of a function of a {Markov} chain and application to
  linkage analysis.
\newblock {\em Advances in Applied Probability}, 31(2):510--531.

\bibitem[Ward et~al., 2022]{Poisson-FOCuS2022}
Ward, K., Dilillo, G., Eckley, I., and Fearnhead, P. (2022).
\newblock {Poisson-FOCuS}: An efficient online method for detecting count
  bursts with application to gamma ray burst detection.
\newblock {\em arXiv preprint arXiv:2208.01494}.

\bibitem[Willsky and Jones, 1976]{willsky1976generalized}
Willsky, A. and Jones, H. (1976).
\newblock A generalized likelihood ratio approach to the detection and
  estimation of jumps in linear systems.
\newblock {\em IEEE Transactions on Automatic control}, 21(1):108--112.

\bibitem[Xie et~al., 2022]{xie2022window}
Xie, L., Moustakides, G.~V., and Xie, Y. (2022).
\newblock Window-limited {CUSUM} for sequential change detection.
\newblock {\em Accepted by IEEE Transactions on Information Theory, arXiv
  preprint arXiv:2206.06777}.

\bibitem[Xie et~al., 2020]{xie2020sequential}
Xie, L., Xie, Y., and Moustakides, G.~V. (2020).
\newblock Sequential subspace change point detection.
\newblock {\em Sequential Analysis}, 39(3):307--335.

\bibitem[Xie et~al., 2021]{xie2021sequential}
Xie, L., Zou, S., Xie, Y., and Veeravalli, V.~V. (2021).
\newblock Sequential (quickest) change detection: Classical results and new
  directions.
\newblock {\em IEEE Journal on Selected Areas in Information Theory},
  2(2):494--514.

\bibitem[Xie et~al., 2012]{xie2012change}
Xie, Y., Huang, J., and Willett, R. (2012).
\newblock Change-point detection for high-dimensional time series with missing
  data.
\newblock {\em IEEE Journal of Selected Topics in Signal Processing},
  7(1):12--27.

\bibitem[Xie and Siegmund, 2013]{xie2013sequential}
Xie, Y. and Siegmund, D. (2013).
\newblock Sequential multi-sensor change-point detection.
\newblock In {\em 2013 Information Theory and Applications Workshop (ITA)},
  pages 1--20. IEEE.

\bibitem[Yakir, 2013]{yakir2013extremes}
Yakir, B. (2013).
\newblock {\em Extremes in random fields: a theory and its applications}.
\newblock John Wiley \& Sons.

\bibitem[Yu et~al., 2020]{yu2020note}
Yu, Y., Padilla, O. H.~M., Wang, D., and Rinaldo, A. (2020).
\newblock A note on online change point detection.
\newblock {\em arXiv preprint arXiv:2006.03283}.

\end{thebibliography}



\end{document}